# HAZARD MODELS WITH VARYING COEFFICIENTS FOR MULTIVARIATE FAILURE TIME DATA[1]

By Jianwen Cai, Jianqing Fan[2], Haibo Zhou[3] and Yong Zhou[4]

*University of North Carolina at Chapel Hill, Princeton University, University of North Carolina at Chapel Hill and Chinese Academy of Science*

Statistical estimation and inference for marginal hazard models with varying coefficients for multivariate failure time data are important subjects in survival analysis. A local pseudo-partial likelihood procedure is proposed for estimating the unknown coefficient functions. A weighted average estimator is also proposed in an attempt to improve the efficiency of the estimator. The consistency and asymptotic normality of the proposed estimators are established and standard error formulas for the estimated coefficients are derived and empirically tested. To reduce the computational burden of the maximum local pseudo-partial likelihood estimator, a simple and useful one-step estimator is proposed. Statistical properties of the one-step estimator are established and simulation studies are conducted to compare the performance of the one-step estimator to that of the maximum local pseudo-partial likelihood estimator. The results show that the one-step estimator can save computational cost without compromising performance both asymptotically and empirically and that an optimal weighted average estimator is more efficient than the maximum local pseudo-partial likelihood estimator. A data set from the Busselton Population Health Surveys is analyzed to illustrate our proposed methodology.

Received May 2004; accepted February 2006.

[1]Supported in part by NIH Grant R01 HL-69720.

[2]Supported in part by NSF Grants DMS-03-55179 and DMS-03-54223 and RGC Grant CUHK4262/01P of HKSAR.

[3]Supported in part by NIH Grant R01 CA-79949.

[4]Supported in part by National Natural Science Foundation of China (NSFC) Grant 10471140.

*AMS 2000 subject classifications.* Primary 62G05; secondary 62N01, 62N02.

*Key words and phrases.* Local pseudo-partial likelihood, marginal hazard model, martingale, multivariate failure time, one-step estimator, varying coefficients.

**1. Introduction.** Multivariate failure time data are encountered in many biomedical studies when related subjects are at risk of a common event or a study subject is at risk of different types of events or recurrence of the same event. Some examples are: epidemiological cohort studies in which the ages of disease occurrence are recorded for members of families; animal experiments where treatments are applied to samples of litter mates; clinical trials in which individual study subjects are followed for the occurrence of multiple events; intervention trials involving group randomization. A common feature of the data in these examples is that the failure times might be correlated. For example, in clinical trials where the patients are followed for repeated recurrent events, the times between recurrences for a given patient may be correlated.

When there is at most one event for each subject and these subjects are mutually independent, the Cox [11] proportional hazards model has commonly been used to assess the effects of the covariates on failure times. For multivariate failure time data, research efforts have been concentrated on marginal hazards models and frailty models. Related literature includes, but is not limited to, Wei, Lin and Weissfeld [28], Lin [20], Cai and Prentice [3, 4] and Spiekerman and Lin [25] for the marginal models and Vaupel, Manton and Stallard [27], Clayton and Cuzick [10], Anderson and Louis [2], Oakes and Jeong [22] and Fan and Li [16] for frailty models. The statistical methods developed for dealing with failure time data typically assume that the covariate effects on the logarithm of the hazard function are linear and that the regression coefficients are constants. These assumptions, however, are primarily made for their mathematical convenience. True associations in practical studies are usually more complex than a simple linear model can capture.

An important extension of the standard regression model with constant coefficient is the varying-coefficient model. The varying-coefficient model addresses an issue frequently encountered by investigators in practical studies. For example, the effect of an exposure variable on the hazard function may change with the level of a confounding covariate. This is traditionally modeled by including an interaction term in the model. Such an approach is a simplification of the true underlying association since a cross product of the exposure and the confounding variable only allows the effect of the exposure to change linearly with the confounding variable. In many studies, however, investigators express the belief that the rate of change is not linear and seek to examine how each level of the exposure interacts with the confounding variables. For example, in a study of cancer risk in uranium miners [23], radon exposure was measured for over 23,000 underground miners in the Czech Republic during 1949–1975. The mining industry's workplace safety measures which affect the inhalation of radon gas, such as ventilation conditions, have changed over the last fifty years. Therefore, the effect of a fixed amount of exposure in the 1950s should not be treated the same as in the 1970s. How to handle this issue is of current active research interest in epidemiology. This leads to a general varying-coefficient model where the coefficient for radon exposure is a function of the calendar year and where this function can be nonlinear over time. Parametric models for the varying-coefficient functions are most efficient if the underlying functions are correctly specified. However, misspecification may cause serious bias and the model constraints may distort the trend in local areas. Nonparametric modeling is appealing in these situations.

Varying-coefficient models have been studied in many non-failure time data settings such as multidimensional nonparametric regression, generalized linear models, analysis of longitudinal data and nonlinear time series. They are particularly appealing in longitudinal studies because they allow one to explore how the effects of the covariates change over time. Related literature includes Hastie and Tibshirani [17], Carroll, Ruppert and Welsh [9] and Cai, Fan and Li [6]. For univariate survival time, the time-varying effect has been carefully studied by Murphy [21], Cai and Sun [7] and Tian, Zuker and Wei [26]. Applications of varying-coefficient models to survival analysis, particularly in the context of multivariate failure time, remain to be studied. New technical challenges arise in dealing with within-cluster dependence and the varying effects of an exposure variable. The local pseudo-partial likelihood in our setting is more sophisticated than that based on the time-varying model. In fact, the latter is no longer a proportional hazards model.

In this paper, we study the marginal hazards model with varying coefficients for multivariate failure time data. The rest of this paper is organized as follows. In Section 2, we formulate the varying-coefficient model and propose local pseudo-partial likelihood procedures for coefficient functions. We also establish asymptotic properties and propose a variance estimator. Further, we consider a computationally efficient one-step procedure and show that it is asymptotically equivalent to the local pseudo-partial likelihood estimator. In Section 3, we propose a weighted average approach to estimate the coefficient functions. We evaluate the proposed procedures through simulation studies and illustrate the proposed approach via an application to the Busselton Population Health Surveys data set in Section 4. Final remarks are made in Section 5. Proofs of theoretical results are given in Section 6.

**2. Marginal hazards model with varying coefficients.** Suppose that there is a random sample of $n$ clusters from an underlying population and that there are $J$ members in each cluster. Let $i$ indicate cluster and $(i,j)$ denote the $j$th member in the $i$th cluster. Let $T_{ij}$ denote the failure time, $C_{ij}$ the censoring time and $X_{ij} = \min(T_{ij}, C_{ij})$ the observed time for member $(i,j)$ ($i=1,\ldots,n$, $j=1,\ldots,J$). Let $\Delta_{ij}$ be an indicator which equals 1 if $X_{ij}$ is a failure time and is 0 otherwise. Varying cluster sizes can be accommodated by defining $T_{ij} \equiv C_{ij} \equiv 0$. Let $\mathcal{F}_{ij}(t)$ represent the failure, censoring and covariate information up to time $t$ for member $(i,j)$ as well as the covariate information for the other members in the $i$th cluster up to time $t$. The marginal hazard function is defined as $\lambda_{ij}(t;\mathcal{F}_{ij}(t)) = \lim_{h\downarrow 0}\frac{1}{h}\mathcal{P}[T_{ij} \leq t+h|T_{ij} > t, \mathcal{F}_{ij}(t)]$. The observed data structure is $\{X_{ij}, \Delta_{ij}, \mathbf{Z}_{ij}(t), V_{ij}(t)\}$ for $i=1,\ldots,n$, where $\mathbf{Z}_{ij}(t) = (Z_{ij1}(t),\ldots,Z_{ijp}(t))^T$ and $V_{ij}(t)$ are two types of covariates, with $V$ being an exposure variable of interest. We assume that the censoring times are independent of the failure times conditional on the covariates and that the observation period is $[0,\tau]$, where $\tau$ is a constant denoting the time for end of the study.

To explore how the effect of the exposure variable $Z$ changes with different levels of a covariate variable $V$, we consider the varying-coefficient model

$$(1) \qquad \lambda_{ij}(t;\mathcal{F}_{ij}) = \lambda_{0j}(t)\exp\{\boldsymbol{\beta}(V_{ij}(t))^T\mathbf{Z}_{ij}(t) + g(V_{ij}(t))\},$$

where $\lambda_{0j}(\cdot)$ is an unspecified baseline hazard function pertaining to the $j$th member of each response vector, $\beta(\cdot)$ is the regression coefficient vector that may be a function of the covariate $V_{ij}$, $g(\cdot)$ is a nonlinear effect of $V_{ij}$ and both $\beta(\cdot)$ and $g(\cdot)$ are unspecified, continuously differentiable functions. Let $N_{ij}(t) = I(X_{ij} \leq t, \Delta_{ij} = 1)$ denote the counting process corresponding to $T_{ij}$ and let $Y_{ij}(t) = I(X_{ij} \geq t)$ denote the 'at risk' indicator process. Set $M_{ij}(t) = N_{ij}(t) - \int_0^t Y_{ij}(s)\lambda_{ij}(s)\,ds$. Note that $M_{ij}(t)$ is a martingale with respect to the marginal filtration $\mathcal{F}_{t,ij} = \sigma\{N_{ij}(s^-), Y_{ij}(s), Z_{ij}(s),\ 0 \leq s \leq t\}$ as well as the union $\sigma$-field $\bigcup_{i=1}^n \mathcal{F}_{t,ij} = \sigma\{N_{ij}(s^-),\ Y_{ij}(s), Z_{ij}(s),\ 0 \leq s \leq t,\ i=1,2,\ldots,n\}$. However, $M_{ij}(t)$ $(i=1,2,\ldots,n,\ j=1,2,\ldots,J)$ is no longer a martingale with respect to the entire union $\sigma$-field $\bigcup_{i=1}^n \bigcup_{j=1}^J \mathcal{F}_{t,ij} = \sigma\{N_{ij}(s^-), Y_{ij}(s), Z_{ij}(s), 0 \leq s \leq t, i=1,2,\ldots,n, j=1,2,\ldots,J\}$ because the observations within a cluster might be mutually dependent.

For ease of presentation, we drop the dependence of covariates on the time $X_{ij}$ with the understanding that the methods and proofs in this paper are applicable to external time-dependent covariates [18]. If all of the observations are independent, then the partial likelihood for model (1) is

$$(2) \quad L(\boldsymbol{\beta}(\cdot),g(\cdot)) = \prod_{j=1}^J \prod_{i=1}^n \left\{\frac{\exp\{\boldsymbol{\beta}(V_{ij})^T\mathbf{Z}_{ij} + g(V_{ij})\}}{\sum_{l \in \mathcal{R}_j(X_{ij})}\exp\{\boldsymbol{\beta}(V_{lj})^T\mathbf{Z}_{lj} + g(V_{lj})\}}\right\}^{\Delta_{ij}},$$

where $\mathcal{R}_j(t) = \{i : X_{ij} \geq t\}$ denotes the set of individuals at risk just prior to time $t$. Since the observations within a cluster are not independent, we refer to (2) as the *pseudo-partial likelihood*. Wei, Lin and Weissfeld [28] considered the parametric counterpart for (2).

2.1. *Local pseudo-partial likelihood estimation.* If the unknown functions $\boldsymbol{\beta}(\cdot)$ and $g(\cdot)$ are parameterized, the parameters can be estimated by maximizing (2). For our nonparametric estimation, the forms of the unknown functions are not available. Directly solving the pseudo-partial likelihood (2) for the unknown functions $\boldsymbol{\beta}(\cdot)$ and $g(\cdot)$ is hardly possible due to the infinite dimension of the unknown parameters. We choose to use the local polynomial method for our nonlinear modeling for the following reasons. First, it is relatively easy to program because existing software on parametric fitting can be modified, via the introduction of a weighting scheme, to deal with local parametric problems. Second, the sampling properties of local polynomial fitting can be derived and efficient semiparametric estimators can be constructed.

Assume that all functions in the components of $\boldsymbol{\beta}(\cdot)$ and $g(\cdot)$ are smooth so that they admit Taylor expansions: for each given $v$ and $u$, where $u$ is close to $v$,

$$\begin{aligned}(3)\quad &\boldsymbol{\beta}(u) \approx \boldsymbol{\beta}(v) + \boldsymbol{\beta}'(v)(u-v) \equiv \boldsymbol{\delta} + \boldsymbol{\eta}(u-v),\\ &g(u) \approx g(v) + g'(v)(u-v) \equiv \alpha + \gamma(u-v),\end{aligned}$$

where $\boldsymbol{\beta}'(u) = d\boldsymbol{\beta}(u)/du$. Substituting these local models into (2), we obtain the following logarithm of the local pseudo-partial likelihood:

$$\begin{aligned}(4)\quad \ell(\gamma, \boldsymbol{\delta}, \boldsymbol{\eta}) = &\sum_{j=1}^{J}\sum_{i=1}^{n} K_h(V_{ij} - v)\\ &\times \Delta_{ij}\bigg\{\boldsymbol{\delta}^T \mathbf{Z}_{ij} + \boldsymbol{\eta}^T \mathbf{Z}_{ij}(V_{ij}-v) + \gamma(V_{ij}-v)\\ &\quad - \log\bigg(\sum_{l \in \mathcal{R}_j(X_{ij})} \exp\{\boldsymbol{\delta}^T \mathbf{Z}_{lj} + \boldsymbol{\eta}^T \mathbf{Z}_{lj}(V_{lj}-v)\\ &\quad\quad\quad\quad\quad\quad\quad\quad + \gamma(V_{lj}-v)\}\\ &\quad\quad\quad\quad\quad\quad\times K_h(V_{lj}-v)\bigg)\bigg\}.\end{aligned}$$

Here $K_h(\cdot) = K(\cdot/h)/h$, $K(\cdot)$ is a probability density called a *kernel function*, and $h$ represents the size of the local neighborhood called a *bandwidth*. The kernel weight is introduced to reflect the fact that the local model (3) is only applied to the data around $v$.

Using counting process notation and letting $\mathbf{X}_{ij}^* = (\mathbf{Z}_{ij}^T, \mathbf{Z}_{ij}^T(V_{ij}-v), V_{ij}-v)^T$ and $\boldsymbol{\xi} = (\boldsymbol{\delta}^T, \boldsymbol{\eta}^T, \gamma)^T$, the local pseudo-partial likelihood function (4) can be expressed as $n \cdot \ell_n(\boldsymbol{\xi}, \infty)$, where

$$\ell_n(\boldsymbol{\xi}, \tau) = n^{-1} \sum_{j=1}^{J} \sum_{i=1}^{n} \int_0^{\tau} K_h(V_{ij} - v)$$
$$\times \left[ \boldsymbol{\xi}^T \mathbf{X}_{ij}^* - \log\left\{ \sum_{l=1}^{n} Y_{lj}(w) \exp(\boldsymbol{\xi}^T \mathbf{X}_{lj}^*) \right.\right.$$
$$\left.\left. \times K_h(V_{lj} - v) \right\} \right] dN_{ij}(w).$$

(5)

Maximizing $\ell(\gamma, \boldsymbol{\delta}, \boldsymbol{\eta})$ in (4) is equivalent to maximizing $\ell_n(\boldsymbol{\xi}, \tau)$ in (5). For a technical reason, following work in the literature, we maximize (5) for a given finite $\tau$.

Let $\hat{\boldsymbol{\xi}}(v) = (\hat{\gamma}(v)^T, \hat{\boldsymbol{\delta}}(v)^T, \hat{\boldsymbol{\eta}}(v)^T)^T$ be the maximizer of (5). Then $\hat{\boldsymbol{\beta}}(v) = \hat{\boldsymbol{\delta}}(v)$ is a local linear estimator for the coefficient function $\boldsymbol{\beta}(\cdot)$ at the point $v$. Similarly, an estimator of $g'(\cdot)$ at the point $v$ is simply the local slope $\hat{\gamma}(v)$, that is, $\hat{g}'(v) = \hat{\gamma}(v)$. The curve $\hat{g}(\cdot)$ can be estimated by integration of the function $\hat{g}'(v)$. Following Hastie and Tibshirani [17], the integration can be approximated by using the trapezoidal rule. The local pseudo-partial likelihood estimator in (4) is particularly easy to compute. It can be implemented by using existing software such as SAS or S-PLUS with the Cox regression procedure. The only difference is that one needs to incorporate the kernel weights in the Cox regression and must repeatedly apply the procedure at a grid of points in the range of the variable $V$.

2.2. *Assumptions and notation.* To express explicitly asymptotic bias and asymptotic variance of the estimator, we introduce some necessary assumptions and notation. Let $\mu_i = \int x^i K(x)\,dx$ and $\nu_i = \int x^i K^2(x)\,dx$ for $i = 0, 1, 2$. Denote $P(w, \mathbf{z}, v) = P(X \geq w | \mathbf{Z} = \mathbf{z}, V = v)$ and $\rho(w, \mathbf{z}, v) = P(w, \mathbf{z}, v) \times \exp\{\boldsymbol{\beta}_0(v)^T \mathbf{z} + g_0(v)\}$. For $k = 0, 1, 2$, define $\mathbf{a}_{jk}(w, v) = f_j(v) E\{\rho(w, \mathbf{Z}_j, v) \times \mathbf{Z}_j^{\otimes k} | V_j = v\}$, where $f_j(\cdot)$ is the density of $V_j$, and $\mathbf{Z}^{\otimes k} = 1$, $\mathbf{Z}$ and $\mathbf{Z}\mathbf{Z}^T$ for $k = 0$, 1 and 2, respectively. Let $\mathbf{a}_{jk}(v) = \int_0^{\tau} \mathbf{a}_{jk}(w, v)\,dw$ and $\mathbf{a}_k(v) = \sum_{j=1}^{J} \mathbf{a}_{jk}(v)$. We will drop the dependence of $\mathbf{a}_k(w, v)$, $\mathbf{a}_{jk}(w, v)$, $\mathbf{a}_{jk}(v)$ and $\mathbf{a}_k(v)$ on $v$ when there is no ambiguity. Finally, let

$$\boldsymbol{\Gamma} = \boldsymbol{\Gamma}(v) = \left\{ \sum_{j=1}^{J} \left( \mathbf{a}_{j2} - \int_0^{\tau} \mathbf{a}_{j1}(w) \mathbf{a}_{j1}(w)^T \mathbf{a}_{j0}(w)^{-1} \lambda_{0j}(w)\,dw \right) \right\}^{-1}$$

and

$$\mathbf{Q} = \begin{pmatrix} \mathbf{Q}_1^{-1} & -(\mathbf{a}_0^{-1})^T \mathbf{Q}_1^{-1} \mathbf{a}_1 \\ -\mathbf{a}_1^T \mathbf{Q}_1^{-1} \mathbf{a}_0^{-1} & (\mathbf{a}_0 - \mathbf{a}_1^T \mathbf{a}_2 \mathbf{a}_1)^{-1} \end{pmatrix},$$

where $\mathbf{Q}_1 = \mathbf{a}_2 - \mathbf{a}_1 \mathbf{a}_1^T \mathbf{a}_0^{-1}$.

Let $\|\cdot\|$ denote the $L_2$-norm and $\|\cdot\|_\Phi$ be the sup-norm of a function or process on a set $\Phi$. The support of the random variable $V$ is denoted by $\mathcal{V}$. For a compact subset $\Phi_V$ of $\mathcal{V}$ and some $\varepsilon > 0$, we define the neighborhood set of $\Phi_{V,\varepsilon}$ as $\Phi_{V,\varepsilon} = \{u : \inf_{v \in \Phi_V} |u - v| \leq \varepsilon\}$. For $k = 0, 1, 2$, let

$$s_{jk}(w, \boldsymbol{\zeta}, v) = f_j(v) \int E[P(w, \mathbf{Z}_j, v)\tilde{\Delta}(y, \boldsymbol{\zeta}_j, w)\mathbf{R}_j^{\otimes k}(y, w)|V_j = v]K(y)\, dy,$$

where $\mathbf{R}_j(y, w) = (\mathbf{Z}_j^T(w), \mathbf{Z}_j^T(w)y, y)^T$ and $\tilde{\Delta}(y, \boldsymbol{\zeta}_j, w) = \exp\{\boldsymbol{\zeta}^T \mathbf{R}_j(y, w) + \boldsymbol{\xi}_0^T \mathbf{R}_j(0, w)\}$, where $\boldsymbol{\xi}_0(\cdot) = (\boldsymbol{\beta}_0^T(\cdot), \boldsymbol{\beta}_0'(\cdot)^T, g_0'(\cdot))^T$.

The following conditions are needed in the proofs of the main results:

(i) The kernel function $K(\cdot) \geq 0$ is a bounded, symmetric function with compact support.

(ii) The functions $\boldsymbol{\beta}(\cdot)$ and $g(\cdot)$ have continuous third derivatives around the point $v$.

(iii) $f_j(\cdot)$ is continuous at the point $v$.

(iv) The conditional probability $P(w, \mathbf{Z}_j(w), \cdot)$ is equicontinuous at $v$ and $\mathbf{Z}_j(w)$ is continuous about $w$ for each $j$.

(v) (a) $nh/\log n \to \infty$ and $nh^5$ is bounded; (b) $\int_0^\tau \lambda_{0j}(t)\, dt < \infty$ for each $j \in \{1, 2, \ldots, J\}$.

(vi) $s_{jr}(t, \boldsymbol{\theta}, v)$, $j = 1, 2, \ldots, J$, $r = 0, 1, 2$, is bounded away from 0 on the product space $[0, \tau] \times \mathbf{C} \times \Phi_{V,\varepsilon}$, that is,

$$\inf_{t \in [0,\tau]} \inf_{\boldsymbol{\theta} \in \mathbf{C}} \inf_{v \in \Phi_{V,\varepsilon}} s_{jr}(t, \boldsymbol{\theta}, v) > 0,$$

where $\boldsymbol{\theta} = (\boldsymbol{\beta}^T, g)$ and

$$\sup_{t \in [0,\tau]} \sup_{\boldsymbol{\theta} \in \mathbf{C}} E\|\mathbf{Z}_j(t)\|^2 \exp(\boldsymbol{\beta}^T \mathbf{Z}_j(t) + g) < \infty$$

for each $j \in \{1, 2, \ldots, J\}$. Meanwhile, $s_{j0}(t, \boldsymbol{\theta}, v)$, $j = 1, 2, \ldots, J$, are continuous functions for $(t, \boldsymbol{\theta}, v) \in [0, \tau] \times \mathbf{C} \times \Phi_{V,\varepsilon}$, uniformly in $t \in [0, \tau]$, and

$$s_{j1}(t, \boldsymbol{\theta}, v) = \frac{\partial}{\partial \boldsymbol{\theta}} s_{j0}(t, \boldsymbol{\theta}, v)$$

and

$$s_{j2}(t, \boldsymbol{\theta}, v) = \frac{\partial^2}{\partial \boldsymbol{\theta}^2} s_{j0}(t, \boldsymbol{\theta}, v).$$

(vii) (Asymptotic variance) The matrix

$$\sum_{j=1}^J \left(\mathbf{a}_{j2} - \int_0^\tau \frac{\mathbf{a}_{j1}(w)\mathbf{a}_{j1}(w)^T}{\mathbf{a}_{j0}(w)}\, d\Lambda_{0j}(w)\right)$$

is positive definite for any $v \in \Phi_{V,\varepsilon}$ and the matrix

$$\mathbf{Q}_2 = \begin{pmatrix} \mathbf{a}_2 & \mathbf{a}_1 \\ \mathbf{a}_1^T & \mathbf{a}_0 \end{pmatrix}$$

is nonsingular at $v \in \Phi_{V,\varepsilon}$.

(viii) The conditional probability $P(u, \mathbf{Z}_j(u), w)$ is equicontinuous in the arguments $(u, w)$ on $[0, \tau] \times \Phi_{W,\varepsilon}$.

(ix) The compact set $\Phi_V \subset \mathcal{W}$ has the following property: $\inf_{u \in \Phi_{W,\varepsilon}} f_j(u) > 0$ for each $j$ and some $\varepsilon > 0$ and $\|f_j\|_{\Phi_V} < \infty$.

(x) The covariate process $\mathbf{Z}_j(t)$ has a continuous sample path in a subset $\mathcal{Z}$ of the continuous function space and $|Z_{ijk}(0)| + \int_0^\tau |dZ_{ijk}(t)| \le B_z$ a.s. for all $i, j, k$ and some constant $B_z < \infty$.

The above conditions will be used for deriving the pointwise convergence properties of $\hat{\boldsymbol{\xi}}$ and demonstrating its asymptotic normality. Conditions (i)–(v) are similar to those in [15] and conditions (vii)–(viii) are similar to conditions C and D of [1]. In order to derive the uniform consistency result, conditions (ix)–(x) are also necessary. From the proofs of the theorems, continuity of $\mathbf{Z}_j(t)$ in assumption (x) can be weakened to $\mathbf{Z}_j(t)$ being left continuous with right-hand limits and $E[\exp\{\beta(V)^T \mathbf{Z}(t)\} \mathbf{Z}(t)^{\otimes k} | V = v]$ and $E(\mathbf{Z}(t)^{\otimes k} | V = v\}$ being continuous functions of $t$ for $k = 0, 1, 2$.

2.3. *Asymptotic properties.* We now establish the asymptotic properties of the local pseudo-partial likelihood estimator. We summarize the results here and provide outlines of the proofs in Section 6. As shown in Section 6, the local pseudo-partial likelihood function $\ell_n(\boldsymbol{\xi}, \tau)$ is concave in $\boldsymbol{\xi}$ and its maximizer exists with probability tending to one. Let $\mathbf{H}$ be a $(2p+1) \times (2p+1)$ diagonal matrix with the first $p$ diagonal elements being 1 and the rest being $h$, where $p$ is the number of elements in $\mathbf{Z}$.

THEOREM 1. *Under conditions* (i)–(viii), *we have*

$$\mathbf{H}\{\hat{\boldsymbol{\xi}}(v) - \boldsymbol{\xi}_0(v)\} \xrightarrow{P} 0,$$

*where* $\boldsymbol{\xi}_0(v) = (\boldsymbol{\beta}_0^T(v), \boldsymbol{\beta}_0'(v)^T, g_0'(v))^T$ *is the vector of the true parameter functions. If, in addition, conditions* (ix)–(x) *are satisfied, then we have the uniform consistency*

$$\sup_{u \in \Phi_V} |\mathbf{H}(\hat{\boldsymbol{\xi}}(u) - \boldsymbol{\xi}_0(u))| \xrightarrow{P} 0,$$

*where $\Phi_V$ is any compact subset of the support of the random variable $V$.*

THEOREM 2. *Assume that conditions* (i)–(viii) *are satisfied. Then the random vector* $(nh)^{-1/2}\{\ell_n'(\boldsymbol{\xi}_0(u), \tau) - \frac{1}{2}h^2\nu_2[(\boldsymbol{\Gamma}^{-1}(u)\boldsymbol{\beta}_0''(u))^T, \mathbf{0}^T, 0]^T\}$ *converges in distribution to a $(2p+1)$-variate normal vector with mean zero and* *covariance matrix* $\mathbf{\Pi}$, *where* $\ell'_n(\boldsymbol{\xi}, \tau) = \partial \ell_n(\boldsymbol{\xi}, \tau)/\partial \boldsymbol{\xi}$, $\mathbf{0}$ *is a p-variate column vector with all entries* 0 *and* $\mathbf{\Pi} = \mathbf{\Pi}_0 + \mathbf{D}$, *in which* $\mathbf{D} = \text{blockdiag}(\mathbf{\Gamma}^{-1}\nu_0, \mathbf{Q}_2\nu_2)$ *and*

$$\mathbf{\Pi}_0 = \sum_{l=1}^{J} \sum_{j=1, j\neq l}^{J} \lim_{n \to \infty} Eh\mathbf{B}_{n1j}(\tau)\mathbf{B}_{n1l}(\tau)^T,$$

*the definitions of* $\mathbf{B}_{n1j}(\tau)$ *and* $\mathbf{B}_{n1l}(\tau)$ *appearing in the proof of Theorem* 2.

THEOREM 3 (Asymptotic normality). *Assume that conditions* (i)–(viii) *are satisfied. Then*

$$\sqrt{nh}\left\{\mathbf{H}(\hat{\boldsymbol{\xi}}(v) - \boldsymbol{\xi}_0(v)) - \frac{1}{2}h^2\mathbf{e}_p\boldsymbol{\xi}''_0(v)\nu_2\right\} \xrightarrow{\mathcal{L}} N(0, \boldsymbol{\Sigma}(\tau, v)),$$

*where* $\mathbf{e}_p$ *is a* $(2p+1) \times (2p+1)$ *matrix with the first* $p \times p$ *elements being* 1 *and the rest being* 0, *and* $\boldsymbol{\Sigma} = \mathbf{A}^{-1}\mathbf{\Pi}(\mathbf{A}^{-1})^T$.

From the expressions for the asymptotic bias and variance matrix $\boldsymbol{\Sigma}$ in Section 6, it can be shown that they can be consistently estimated by

(6) $\quad \widehat{\mathbf{A}}_n^{-1}(\tau, v)\widehat{\mathbf{B}}_n(\tau, v) \quad \text{and} \quad (nh)^{-1}\widehat{\mathbf{A}}_n^{-1}(\tau, v)\widehat{\mathbf{\Pi}}_n(\tau, v)\widehat{\mathbf{A}}_n^{-1}(\tau, v),$

where

$$\widehat{\mathbf{A}}_n(\tau, v) = \frac{1}{n}\sum_{i=1}^{n}\sum_{j=1}^{J}\int_0^\tau K_h(V_{ij} - v)\left(\frac{\widehat{S}_{nj2}(w, v)}{\widehat{S}_{nj0}(w, v)} - \widehat{E}_j(w, v)^{\otimes 2}\right)dN_{ij}(w),$$

$$\widehat{\mathbf{B}}_n(\tau, v) = \frac{1}{nh}\sum_{i=1}^{n}\sum_{j=1}^{J}\int_0^\tau K_h(V_{ij} - v)(\mathbf{U}^*_{ij}(w) - \widehat{E}_j(w, v))\,dN_{ij}(w),$$

$$\widehat{\mathbf{\Pi}}_n(\tau, v) = \frac{1}{nh}\sum_{i=1}^{n}\left\{\sum_{j=1}^{J}\int_0^\tau K_h(V_{ij} - v)(\mathbf{U}^*_{ij}(w) - \widehat{E}_j(w, v))\,d\widehat{M}_{ij}(w)\right\}^{\otimes 2}$$

with $\widehat{S}_{njk}(w, v) = \frac{1}{n}\sum_{i=1}^{n} K_h(V_{ij} - v)Y_{ij}(w)\exp(\hat{\boldsymbol{\xi}}^T(v)\mathbf{X}^*_{ij}(w))(\mathbf{U}^*_{ij}(w))^{\otimes k}$ for $k = 0, 1, 2$, $\mathbf{U}^*_{ij} = \mathbf{H}^{-1}\mathbf{X}^*_{ij}$, $\widehat{E}_j(w, v) = \widehat{S}_{nj1}(w, v)/\widehat{S}_{nj0}(w, v)$ and $\widehat{M}_{ij}(t) = N_{ij}(t) - \int_0^t \hat{\lambda}_{ij}(s)ds$, in which $\hat{\lambda}_{ij}(s) = \hat{\lambda}_{0j}(s)\exp\{\hat{\boldsymbol{\beta}}(V_{ij}(s))Z_{ij}(s) + \hat{g}(V_{ij}(s))\}$ and $\hat{\lambda}_{0j}(s)$ is given in the following section.

2.4. *Estimation of the baseline hazard function.* With estimators of $\boldsymbol{\beta}(\cdot)$ and $g(\cdot)$, we can estimate the baseline hazard function by using a kernel smoothing,

$$\hat{\lambda}_{0j}(t) = \int W_b(t - x)\,d\widehat{\Lambda}_{0j}(x),$$

where $W_b$ is a given kernel function and $b$ is a given bandwidth. The cumulative hazard function $\Lambda_{0j}(\cdot)$ can be estimated by

$$\widehat{\Lambda}_{0j}(t) = \frac{1}{n}\sum_{i=1}^{n}\int_0^t \frac{dN_{ij}(w)}{n^{-1}\sum_{l=1}^{n} Y_{lj}(w)\exp(\hat{\boldsymbol{\beta}}(V_{lj})^T \mathbf{Z}_{lj}(w) + \hat{g}(V_{lj}))}.$$

The properties of $\widehat{\Lambda}_{0j}(\cdot)$ and $\hat{\lambda}_{0j}(\cdot)$ are summarized in the following theorem and an outline of the proof is provided in Section 6.

THEOREM 4. *Under conditions* (i)–(x), *we have*
$$\widehat{\Lambda}_{0j}(t) \longrightarrow \Lambda_{0j}(t) \quad and \quad \hat{\lambda}_{0j}(t) \longrightarrow \lambda_{0j}(t),$$
*uniformly on* $(0,\tau]$ *in probability.*

To investigate the asymptotic properties of the estimated cumulative hazard function, we assume, for simplicity, that $g(V) = 0$. The function $g(\cdot)$ needs to be estimated by integrating its derivative estimator from the partial likelihood; hence, its asymptotic properties are challenging to obtain. When $g(\cdot) = 0$ our task is somewhat simplified. The generalized Breslow estimator for $\Lambda_{0j}(t)$ is given by

$$\widehat{\Lambda}_{0j}(t) = \frac{1}{n}\sum_{i=1}^{n}\int_0^t \frac{dN_{ij}(w)}{n^{-1}\sum_{l=1}^{n} Y_{lj}(w)\exp\{\hat{\boldsymbol{\beta}}(V_{lj})^T \mathbf{Z}_{lj}(w)\}}.$$

Write $W_j(t) = n^{1/2}\{\widehat{\Lambda}_{0j}(t) - \Lambda_{0j}(t)\}$. This is a stochastic process defined in the metric space $\Omega = C[0,\tau]$ with the norm $\rho(f,g) = \max_{1\leq j \leq J} \sup_{0\leq t \leq \tau} |f_j(t) - g_j(t)|$.

THEOREM 5. *Assume that conditions* (i)–(x) *are satisfied and let* $nh^4 \to 0$. *Then the random process vector* $\mathbf{W}(t) = (W_1(t),\ldots,W_J(t))$ *converges weakly to a zero-mean Gaussian random field* $\mathbf{G}(t)$.

REMARK 1. The covariance structure of the Gaussian field $\mathbf{G}(t)$ is very complex. It is very difficult to directly calculate this covariance by asymptotic methods. The wild bootstrap provides a useful method for computing the covariance or approximating the distribution of $\mathbf{G}(t)$ (see [25]).

REMARK 2. Theorem 5 shows that the estimator $\widehat{\Lambda}_{0j}(t)$ is root-$n$ consistent if the nonparametric estimators are undersmoothed. This means that in practical applications, one uses the right amount of smoothing for estimating coefficient functions and then chooses a smaller amount of smoothing for estimating the cumulative hazard functions. The situation here is very different from the one-step likelihood estimation of Carroll et al. [8], but similar to their one-step procedure.

2.5. *One-step local pseudo-partial likelihood estimator.* To estimate the functions $\boldsymbol{\beta}(\cdot)$ and $g(\cdot)$ over an interval of interest, we usually need to maximize the local pseudo-partial likelihood (5) at hundreds of points. This can be very computationally intensive. In addition, for certain given $v$, the local pseudo-partial likelihood estimator might not exist, due to a limited amount of data around $v$. These drawbacks make computing the local pseudo-partial likelihood estimator over an interval less appealing. We consider the following one-step estimator as a feasible alternative.

To facilitate notation, we drop the dependence of $\ell_n(\boldsymbol{\xi}, \tau)$ on $\tau$. The local pseudo-partial likelihood estimator $\hat{\boldsymbol{\xi}}$ satisfies $\ell'_n(\boldsymbol{\xi}) = 0$. For a given initial estimator $\hat{\boldsymbol{\xi}}_0$, by Taylor expansion, we have

$$\ell'_n(\hat{\boldsymbol{\xi}}_0) + \ell''_n(\hat{\boldsymbol{\xi}}_0)(\hat{\boldsymbol{\xi}} - \hat{\boldsymbol{\xi}}_0) \approx 0.$$

Thus, the one-step estimator $\hat{\boldsymbol{\xi}}_{\text{os}}$ is defined as

(7) $$\hat{\boldsymbol{\xi}}_{\text{os}} = \hat{\boldsymbol{\xi}}_0 - \{\ell''_n(\hat{\boldsymbol{\xi}}_0)\}^{-1} \ell'_n(\hat{\boldsymbol{\xi}}_0).$$

In the Newton–Raphson algorithm, the above equation is iterated until convergence is achieved. As shown in Section 6, the function $\ell_n(\boldsymbol{\xi})$ is concave. Hence, its maximizer exists and is unique when $\ell_n(\boldsymbol{\xi})$ is strictly concave. In practice, we do not have to iterate (7) until convergence is achieved—once, or a few times, will suffice. Robinson [24] gives results on the distance between the estimators based on a few iteration steps and the maximum likelihood estimator. A natural question arises as to how good the initial estimator $\hat{\boldsymbol{\xi}}_0$ has to be in order for the one-step estimator to have the same performance as the maximum local pseudo-partial likelihood estimator. It is not hard to show that a sufficient condition is

(8) $$\mathbf{H}(\hat{\boldsymbol{\xi}}_0 - \boldsymbol{\xi}_0) = O_P(h^2 + (nh)^{-1/2});$$

see Fan and Chen [13] for a derivation in the local likelihood context. When condition (8) is not satisfied, a multiple-step estimator is needed. By repeatedly applying the one-step result $k$ times, as in [24], condition (8) can be relaxed to $\mathbf{H}(\hat{\boldsymbol{\xi}}_0 - \boldsymbol{\xi}_0) = O_P\{(h^2 + (nh)^{-1/2})^{1/k}\}$.

Cai, Fan and Li [6] provide a useful strategy on the choice of initial estimators in the context of generalized linear models and their idea can be adapted to the current setting. The idea is to exploit the smoothness of nonparametric functions. The strategy is as follows. Compute the local pseudo-partial likelihood estimates at a few fixed points. Use these estimates as the initial values of their nearest grid points and obtain the one-step estimates at these grid points. Use the newly computed one-step estimates as the initial values of their nearest grid points to compute the one-step estimates and propagate until the one-step estimates at all grid points are computed. For example, in our simulation studies, we evaluate the functions at $n_{\text{grid}} = 200$ grid points and are willing to compute the maximum local pseudo-partial likelihood at five distinct points. A sensible placement of these points is $w_{20}, w_{60}, w_{100}, w_{140}$ and $w_{180}$. We shall use, for instance, $\hat{\boldsymbol{\beta}}(w_{60})$ as an initial value for calculating the one-step estimates for $\hat{\boldsymbol{\beta}}(w_{59})$ and $\hat{\boldsymbol{\beta}}(w_{61})$ and then proceed to use the resulting estimates as the initial values for calculating the one-step estimates for $\hat{\boldsymbol{\beta}}(w_{58})$ and $\hat{\boldsymbol{\beta}}(w_{62})$, respectively. We continue this process until all the one-step estimates at $w_i$ for $i = 40, \ldots, 79$ are calculated.

**3. Weighted average estimator.** An alternative approach is to fit a varying-coefficient model for each failure type, that is, for event type $j$, fitting the model

$$\lambda_{ij}(t; \mathcal{F}_{ij}) = \lambda_{0j}(t) \exp\{\boldsymbol{\beta}_j(V_{ij}(t))^T \mathbf{Z}_{ij}(t) + g_j(V_{ij}(t))\}, \qquad \text{for } i = 1, \ldots, n,$$

resulting in $\hat{\boldsymbol{\xi}}_j(v)$ for estimating $\boldsymbol{\xi}_j(v) = (\boldsymbol{\beta}_j^T(v), (\boldsymbol{\beta}_j'(v))^T, g_j'(v))$. Under model (1), we have $\boldsymbol{\xi}_1 = \boldsymbol{\xi}_2 = \cdots = \boldsymbol{\xi}_J = \boldsymbol{\xi}$. Thus, we can estimate $\boldsymbol{\xi}(v)$ by a linear combination $c_1 \hat{\boldsymbol{\xi}}_1(v) + \cdots + c_J \hat{\boldsymbol{\xi}}_J(v)$ with $\sum_{j=1}^J c_j = 1$. Weights $\mathbf{c}_j$ can be chosen to optimize the performance. Note that the weights $c_j$ can be generalized to a matrix $\mathbf{C}_j$ to allow for different linear combinations for different components of $\boldsymbol{\xi}(v)$, that is, the linear combination can be generalized to $\mathbf{C}_1 \hat{\boldsymbol{\zeta}}_1 + \cdots + \mathbf{C}_J \hat{\boldsymbol{\zeta}}_J$ with $\mathbf{C}_1 + \cdots + \mathbf{C}_J = \text{diag}(1, \ldots, 1)$ being the identity matrix.

In order to establish the asymptotic distribution of the weighted average estimator, we need to derive the asymptotic distribution of $\hat{\Psi}(v) = (\hat{\boldsymbol{\xi}}_1^T, \ldots, \hat{\boldsymbol{\xi}}_J^T)^T$. We define $\Psi(v)$ and $\Psi''(v)$ similarly to $\hat{\Psi}(v)$, except that the $\hat{\boldsymbol{\xi}}_j$ are replaced by $\boldsymbol{\xi}_j$ and $\boldsymbol{\xi}_j''$, respectively, for $j = 1, 2, \ldots, J$. Using arguments similar to those used for Theorems 2 and 3, it can be shown that the following theorem holds.

THEOREM 6. *Under the conditions of Theorem 2, we have that*

$$\sqrt{nh} \left\{ \mathbf{I}_p \otimes \mathbf{H}[\widehat{\Psi}(v) - \Psi(v)] - \frac{h^2 \nu_2}{2} \mathbf{R} \Psi''(v) \right\}$$

*is asymptotically normal with mean* 0 *and covariance matrix* $\boldsymbol{\Sigma}^* = (G_{kl}(\boldsymbol{\xi}_k, \boldsymbol{\xi}_l))$ *for* $k, l = 1, \ldots, J$, *where* $\mathbf{R} = \text{diag}(\mathbf{R}_1, \ldots, \mathbf{R}_J)$ *and* $\mathbf{R}_j$ *is a* $(2p+1) \times (2p+1)$ *matrix with the first* $p \times p$ *elements being* $\mathbf{I}_p$ *(an identity matrix) and rest equal to* 0. $G_{kl}(\boldsymbol{\xi}_k, \boldsymbol{\xi}_l)$ *is defined at the end of this section.*

The asymptotic normality of the weighted average estimator follows easily from Theorem 6. For example, suppose that we are interested in estimating the $k$th entry of $\boldsymbol{\beta}$. Write $\mathbf{1}_k$ as a $(2p+1)$-variate vector with the $k$th entry being 1 and the rest being 0 and let $\mathcal{C} = \operatorname{diag}(\mathbf{1}_k^T,\ldots,\mathbf{1}_k^T) = \mathbf{I}_J \otimes \mathbf{1}_k^T$. Then it follows from Theorem 6 that

$$(nh)^{1/2}\left\{(\hat{\boldsymbol{\beta}}_w(v) - \boldsymbol{\beta}_0(v)) - \frac{h^2\nu_2}{2}\mathcal{C}\mathbf{R}\Psi''(v)\right\} \xrightarrow{\mathcal{L}} N(0,\boldsymbol{\Sigma}_w),$$

where $\hat{\boldsymbol{\beta}}_w = (\hat{\boldsymbol{\beta}}_{k1},\ldots,\hat{\boldsymbol{\beta}}_{kJ})^T$, $\boldsymbol{\beta}_0 = (\boldsymbol{\beta}_{k1},\ldots,\boldsymbol{\beta}_{kJ})^T$, and $\boldsymbol{\Sigma}_w = \mathcal{C}^T\boldsymbol{\Sigma}^*\mathcal{C}$, $\hat{\boldsymbol{\beta}}_{kj}$ and $\boldsymbol{\beta}_{kj}$ being the $k$th entry of $\hat{\boldsymbol{\beta}}_j$ and $\boldsymbol{\beta}_j$, respectively. The optimal weight which minimizes $\mathbf{c}^T\boldsymbol{\Sigma}_w\mathbf{c}$ with $\sum_{j=1}^J c_j = 1$ is $\mathbf{c}_k = (\mathbf{e}^T\boldsymbol{\Sigma}_w^{-1}\mathbf{e})^{-1}\boldsymbol{\Sigma}_w^{-1}\mathbf{e}$. Since the failure times for different types of failures are usually mutually dependent, $\hat{\boldsymbol{\xi}}_j$ $(j=1,\ldots,J)$ are likely to be dependent; hence, the variance $\boldsymbol{\Sigma}_w$ is not necessarily diagonal. This implies that the optimal weight is unlikely to be uniform. In other words, the weighted average estimator with the optimal weight is generally more efficient than the estimator with the "working independence" weight. This is supported by the simulation results displayed in Table 2.

We now give the expressions for the asymptotic variance-covariance matrix from Theorem 6 and its estimate. From Theorem 3, it is easy to show that the asymptotic covariance matrix between $(nh)^{1/2}\mathbf{H}(\hat{\boldsymbol{\xi}}_k(v) - \boldsymbol{\xi}_{k0}(v))$ and $(nh)^{1/2}\mathbf{H}(\hat{\boldsymbol{\xi}}_l(v) - \boldsymbol{\xi}_{l0}(v))$ is given by

$$G_{kl}(\boldsymbol{\xi}_k,\boldsymbol{\xi}_l) = \mathbf{A}_k^{-1}(\boldsymbol{\xi}_k)\lim_{n\to\infty}E\{\Pi_{1k}(\boldsymbol{\xi}_k)\Pi_{1l}(\boldsymbol{\xi}_l)\}\mathbf{A}_l^{-1}(\boldsymbol{\xi}_l),$$

where $\Pi_{jk}(\boldsymbol{\xi}_k) = \int_0^\tau K_h(V_{jk}-v)[U_{jk}^* - s_{k1}(w,\boldsymbol{\zeta},v)/s_{k0}(w,\boldsymbol{\zeta},v)]\,dM_{jk}(w)$, $\boldsymbol{\zeta} = \mathbf{H}(\boldsymbol{\xi} - \boldsymbol{\xi}_0)$ and $s_{kd}(w,\boldsymbol{\zeta},v)$, $d = 0,1$, are defined as in Section 6. From the definition of $\Pi_{jk}(\boldsymbol{\xi}_k)$, it is natural to estimate $\lim_{n\to\infty}E\{\Pi_{k1}(\boldsymbol{\xi}_k)\Pi_{l1}(\boldsymbol{\xi}_l)\}$ by

$$(9) \qquad \widehat{D}_{kl}(\hat{\boldsymbol{\xi}}_k,\hat{\boldsymbol{\xi}}_l) = n^{-1}\sum_{j=1}^n W_{jk}(\hat{\boldsymbol{\xi}}_k)W_{jl}^T(\hat{\boldsymbol{\xi}}_l),$$

where

$$W_{jk}(\boldsymbol{\xi}_k) = \Delta_{jk}\left\{U_{jk}^*(X_{jk}) - \frac{\widehat{S}_{nj1}(X_{jk},v)}{\widehat{S}_{nj0}(X_{jk},v)}\right\}K_h(V_{jk}-v)$$

$$(10) \qquad -\sum_{m=1}^n \frac{\Delta_{mk}Y_{jk}(X_{mk})\exp\{\hat{\beta}_k^T(V_{jk}(X_{mk}))Z_{jk}(X_{mk}) + \hat{g}_k(V_{mk}(X_{mk}))\}}{\sum_{i=1}^n Y_{ik}(X_{mk})\exp\{\hat{\beta}_k^T(V_{ik}(X_{mk}))Z_{ik}(X_{mk}) + \hat{g}_k(V_{mk}(X_{mk}))\}}$$

$$\times \left\{U_{jk}^*(X_{mk}) - \frac{\widehat{S}_{nj1}(X_{mk},v)}{\widehat{S}_{nj0}(X_{mk},v)}\right\}K_h(V_{jk}-v).$$

Obviously, $\mathbf{A}_j(\boldsymbol{\xi})$ can be estimated by

$$\widehat{\mathbf{A}}_j(\boldsymbol{\xi}) = \frac{1}{n}\sum_{i=1}^n \int_0^\tau K_h(V_{ij}-v)\frac{\widehat{S}_{nj2}(w,v)\widehat{S}_{nj0}(w,v) - \widehat{S}_{nj1}^{\otimes 2}(w,v)}{(\widehat{S}_{nj0}(w,v))^2}\,dN_{ij}(w).$$

Write
$$\widehat{G}_{kl}(\boldsymbol{\xi}_k, \boldsymbol{\xi}_l) = \widehat{\mathbf{A}}_k^{-1}(\boldsymbol{\xi}_k)\widehat{D}_{kl}(\hat{\boldsymbol{\xi}}_k, \hat{\boldsymbol{\xi}}_l)\widehat{\mathbf{A}}_l^{-1}(\boldsymbol{\xi}_l).$$

By means of some tedious proofs, we can show that $\widehat{G}_{kl}(\boldsymbol{\xi}_k, \boldsymbol{\xi}_l)$ is a consistent estimator of $G_{kl}(\boldsymbol{\xi}_k, \boldsymbol{\xi}_l)$. Hence, the covariance matrix of $(\hat{\boldsymbol{\xi}}_1, \ldots, \hat{\boldsymbol{\xi}}_J)$ can be consistently estimated by $\widehat{\boldsymbol{\Sigma}}^* = (nh)^{-1}(\widehat{G}_{ij}(\hat{\boldsymbol{\xi}}_i, \hat{\boldsymbol{\xi}}_j))_{i,j=1}^J$. These results provide a basis for simultaneous inferences about the $\boldsymbol{\xi}_j$, $j = 1, 2, \ldots, J$, as well as for the weighted average estimator $\sum_{j=1}^J c_j \hat{\boldsymbol{\xi}}_j$ for $\boldsymbol{\xi}$.

## 4. Numerical examples.

4.1. *Simulations.* We perform a series of simulation studies to evaluate the performance of the proposed estimation method. Multivariate failure times are generated from a multivariate extension of the model of Clayton and Cuzick [10] in which the joint survival function of $(T_1, \ldots, T_J)$ given $(\mathbf{Z}_1, \ldots, \mathbf{Z}_J)$ and $(V_1, \ldots, V_J)$ is

$$(11) \quad F(t_1, \ldots, t_J; \mathbf{Z}_1, \ldots, \mathbf{Z}_J, V_1, \ldots, V_J) = \left\{\sum_{j=1}^J S_j(t_j)^{-\theta} - (J-1)\right\}^{-1/\theta},$$

where $J$ takes integer values and $S_j(t)$ is the marginal survival probability for the $j$th member, depending on covariates $\mathbf{Z}_j$ and $V_j$. Note that $\theta$ is a parameter which represents the degree of dependence of $T_i$ and $T_j$, $i, j = 1, 2, \ldots, J$. The relationship between Kendall's $\tau$ and $\theta$ is $\tau = \theta/(2 + \theta)$. Specifically, $\theta = 0.25$ and $\theta = 4$ represent weak and strong positive dependence, respectively, with $\theta \to 0$ giving independence and $\theta \to \infty$ giving maximal positive dependence. In our simulation, $\theta$ was chosen to be 0.25, 1.5 and 4.0, which correspond to low, moderate and high positive dependence, respectively. The Gaussian kernel function is used for the estimates.

In our first set of simulations, we examine the performance of the local pseudo-partial likelihood estimators. We consider the marginal distribution of $T_{ij}$ to be exponential with failure rate

$$(12) \quad \lambda_{ij}(t) = \lambda_{0j}(t)\exp\{\beta(V_{ij})\mathbf{Z}_{ij} + g(V_{ij})\}.$$

We choose the baseline hazard function to depend on time, specifically, $\lambda_{0j}(t) = 4t^3 \lambda_{0j}^*$, where $j = 1, 2, 3$. We take $\lambda_{0j}^*$ to be 0.2, 1.0 and 1.5 for $j = 1, 2$ and 3, respectively. Failure times $(t_{i1}, t_{i2}, t_{i3})$, $i = 1, 2, \ldots, n$, are generated from the distribution function (11) with marginal distribution (12). We generate the covariate vector $\mathbf{Z}_{ij} = (Z_{ij1}, Z_{ij2}, \ldots, Z_{ijp})^T$ from a multivariate normal distribution with marginal mean 0, standard deviation 5 and correlation between $Z_{ijl}$ and $Z_{ijk}$ equal to $\rho^{l-k}$, where $\rho = 1/\sqrt{5}$. We consider $p = 2$ and the varying coefficients

$$\beta_1(V) = 0.5V(1.5 - V), \qquad \beta_2(V) = \sin(2V)$$

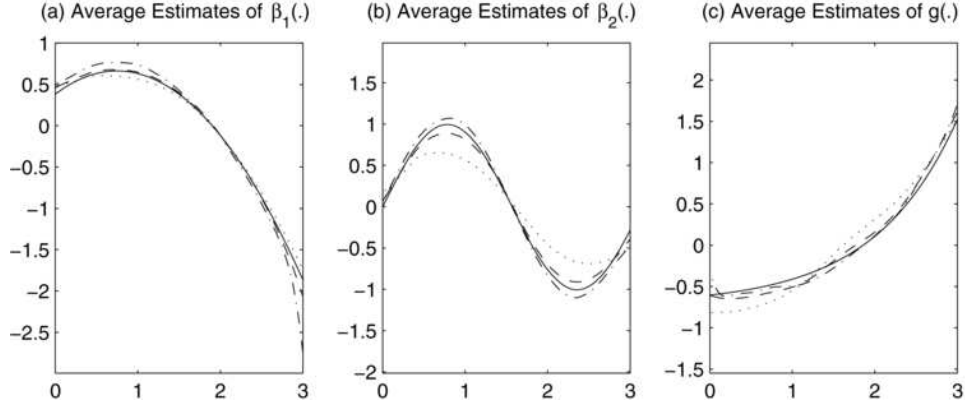

FIG. 1. *Simulation results for the local pseudo-partial likelihood estimator with* 200 *clusters and the dependence parameter of failure time* $\theta = 0.25$. (a), (b) *and* (c) *provide the average estimates of* $\beta_1(\cdot)$, $\beta_2(\cdot)$ *and* $g(\cdot)$ *for heavily censored data with* $c = 2$, *respectively. Solid curves are true functions. Three bandwidths are used: dash-dotted curves for bandwidth* $h = 0.1$, *dashed curves for* $h = 0.2$ *and dotted curves for* $h = 0.4$.

and

$$g(V) = 0.5\{e^{V-1.5} - e^{-1.5}\},$$

where $V$ is generated from a uniform distribution over $[0, 3]$. In our simulation, using a similar derivation as in [5], with given covariates $(\mathbf{Z}_{ij}, V_{ij})$, $j = 1, 2, 3$, the failure times $(t_{i1}, t_{i2}, t_{i3})$ are generated from independent uniform random variables $(w_{i1}, w_{i2}, w_{i3})$ as follows:

$$t_{i1} = [-\log(1 - w_{i1})\Upsilon(V_{i1}, \mathbf{Z}_{i1}, \lambda_{01}^*)]^{1/4},$$
$$t_{i2} = [\theta \log\{1 - a_{i1} + a_{i1}(1 - w_{i2})^{(\theta^{-1}+1)^{-1}}\}\Upsilon(V_{i2}, \mathbf{Z}_{i2}, \lambda_{02}^*)]^{1/4},$$
$$t_{i3} = [\theta \log\{1 - (a_{i1} + a_{i2}) + (a_{i1} + a_{i2} - 1)(1 - w_{i3})^{(\theta^{-1}+2)^{-1}}\}$$
$$\times \Upsilon(V_{i3}, \mathbf{Z}_{i3}, \lambda_{03}^*)]^{1/4}.$$

Here $a_{il} = (1 - w_{il})^{-\theta}$ for $l = 1, 2$ and $i = 1, 2, \ldots, n$ and $\Upsilon(V, \mathbf{Z}, \lambda^*) = \exp\{\boldsymbol{\beta}(V)\mathbf{Z} + g(V)\}/\lambda^*$. Censoring times $C_{ij}$ are generated from the uniform distribution over $(0, c)$, where $c$ is a constant which is set to control the censoring rate. There is approximately 10% censoring when $c = 5$ and approximately 30% censoring when $c = 2$. For each of the configurations studied, 500 simulations were carried out.

Table 1 summarizes the simulation results for the local pseudo-partial likelihood estimator of $\boldsymbol{\beta}(\cdot)$ and $g'(\cdot)$ with the number of clusters being 200, $\theta = 0.25$ and $c = 2$. We present the estimates of the functions evaluated at $v = 0.5, 1.0, 1.5, 2.0$ and $2.5$. The bandwidths we considered were $h = 0.075$, 0.1, 0.15, 0.2 and 0.4. The averages of the 500 estimates for $\beta_1(v)$, $\beta_2(v)$ and $g'(v)$ subtracting their true values are given in the "bias" columns and the standard deviations of the 500 estimates are given in the corresponding SD columns. The SE columns give the averages of the estimated standard errors. Figure 1 provides the average estimates for $\beta_1(\cdot)$, $\beta_2(\cdot)$ and $g(\cdot)$ based on different bandwidths. It gives us an idea of how large the biases are for different bandwidths. From Table 1, we can also see that as the bandwidth increases, the variance decreases. As expected, with a large bandwidth $h = 0.4$, the bias is large and the variance is small. Note that the absolute biases exhibit a U-shape in Table 1. This is unusual, but can happen. The bias depends on function values in a local neighborhood and is continuous in $h$.

TABLE 1
*Summary of simulation results based on local pseudo-partial likelihood procedures*

| $v$ | $h$ | $\hat{\beta}_1(\cdot)$ bias | SE | SD | $\hat{\beta}_2(\cdot)$ bias | SE | SD | $\hat{g}'(\cdot)$ bias | SE | SD |
|---|---|---|---|---|---|---|---|---|---|---|
| 0.5 | 0.075 | −0.071 | 0.194 | 0.262 | −0.179 | 0.266 | 0.370 | −0.059 | 1.931 | 1.501 |
|  | 0.100 | −0.036 | 0.158 | 0.177 | −0.095 | 0.215 | 0.266 | −0.011 | 1.118 | 0.956 |
|  | 0.150 | −0.007 | 0.121 | 0.133 | −0.004 | 0.160 | 0.175 | 0.003 | 0.538 | 0.493 |
|  | 0.200 | 0.028 | 0.096 | 0.101 | 0.077 | 0.121 | 0.125 | 0.044 | 0.328 | 0.295 |
|  | 0.400 | 0.077 | 0.076 | 0.086 | 0.252 | 0.087 | 0.108 | 0.247 | 0.167 | 0.123 |
| 1.0 | 0.075 | −0.041 | 0.187 | 0.252 | −0.173 | 0.279 | 0.354 | 0.085 | 1.817 | 1.495 |
|  | 0.100 | −0.009 | 0.153 | 0.188 | −0.094 | 0.218 | 0.259 | 0.079 | 1.140 | 0.933 |
|  | 0.150 | 0.004 | 0.115 | 0.118 | 0.006 | 0.156 | 0.164 | −0.007 | 0.456 | 0.451 |
|  | 0.200 | 0.027 | 0.092 | 0.093 | 0.101 | 0.113 | 0.124 | −0.066 | 0.305 | 0.257 |
|  | 0.400 | 0.108 | 0.063 | 0.064 | 0.382 | 0.066 | 0.082 | −0.055 | 0.132 | 0.093 |
| 1.5 | 0.075 | −0.004 | 0.174 | 0.231 | −0.031 | 0.181 | 0.236 | 0.066 | 1.843 | 1.527 |
|  | 0.100 | 0.016 | 0.142 | 0.164 | −0.015 | 0.145 | 0.157 | 0.032 | 1.040 | 0.949 |
|  | 0.150 | 0.019 | 0.110 | 0.114 | −0.007 | 0.116 | 0.115 | 0.035 | 0.533 | 0.496 |
|  | 0.200 | 0.023 | 0.092 | 0.094 | 0.010 | 0.096 | 0.097 | 0.035 | 0.329 | 0.304 |
|  | 0.400 | 0.075 | 0.061 | 0.060 | 0.063 | 0.065 | 0.068 | −0.181 | 0.131 | 0.095 |
| 2.0 | 0.075 | 0.097 | 0.196 | 0.277 | 0.145 | 0.235 | 0.326 | 0.140 | 1.886 | 1.503 |
|  | 0.100 | 0.076 | 0.167 | 0.195 | 0.070 | 0.192 | 0.224 | 0.104 | 1.075 | 0.956 |
|  | 0.150 | 0.047 | 0.129 | 0.142 | 0.004 | 0.143 | 0.139 | 0.078 | 0.566 | 0.521 |
|  | 0.200 | 0.037 | 0.107 | 0.112 | −0.044 | 0.114 | 0.117 | 0.091 | 0.355 | 0.339 |
|  | 0.400 | 0.049 | 0.072 | 0.076 | −0.254 | 0.068 | 0.076 | −0.098 | 0.138 | 0.116 |
| 2.5 | 0.075 | 0.206 | 0.315 | 0.413 | 0.141 | 0.261 | 0.329 | 0.470 | 1.982 | 1.635 |
|  | 0.100 | 0.135 | 0.260 | 0.305 | 0.057 | 0.206 | 0.251 | 0.210 | 1.160 | 1.058 |
|  | 0.150 | 0.074 | 0.200 | 0.216 | −0.004 | 0.151 | 0.166 | 0.095 | 0.633 | 0.570 |
|  | 0.200 | 0.039 | 0.166 | 0.179 | −0.074 | 0.119 | 0.128 | 0.005 | 0.385 | 0.384 |
|  | 0.400 | −0.014 | 0.125 | 0.130 | −0.222 | 0.087 | 0.098 | −0.261 | 0.189 | 0.181 |

If the bias associated with $h = 0.075$ is negative and the bias associated with $h = 0.4$ is positive, as in this example, then the bias necessarily crosses zero and the U-shape absolute biases emerge. Table 1 also shows that when the bandwidth is 0.15 and above, the average SE is close to SD, which indicates good performance of the variance estimator. We have also examined the situations involving moderate and high dependence of the failure times, using $\theta = 1.5$ and 4, respectively, as well as lighter censoring (with $c = 5$). The conclusions are similar.

We also examined the performance of the weighted average estimator. The results are presented in Table 2. The rows labeled "W" are those based on the weighted average estimates for $\beta_p$ with optimal weight $\mathbf{c}_p = (\mathbf{e}^T \widehat{\boldsymbol{\Sigma}}_p^{-1} \mathbf{e})^{-1} \widehat{\boldsymbol{\Sigma}}_p^{-1} \mathbf{e}$, where $\widehat{\boldsymbol{\Sigma}}_p$ is the estimator of the asymptotic variance-covariance matrix of $(\hat{\beta}_{p1}, \hat{\beta}_{p2}, \hat{\beta}_{p3})^T$. The maximum local pseudo-partial likelihood estimates are indicated by "P". The performance of an estimator $\hat{\boldsymbol{\beta}}(\cdot)$ is assessed via the square root of average square errors (RASE),

$$(13) \qquad \text{RASE} = \left( \frac{1}{n_{\text{grid}}} \sum_{k=1}^{n_{\text{grid}}} [\hat{\beta}(w_k) - \beta(w_k)]^2 \right)^{1/2},$$

where $\{w_k, k = 1, \ldots, n_{\text{grid}}\}$ are the grid points at which the functions $\boldsymbol{\beta}(\cdot)$ are estimated. In our simulations, we take $n_{\text{grid}} = 200$.

For the local pseudo-partial likelihood estimates, we used bandwidth $h = 0.15$. For the weighted average estimates, a bandwidth of $h = 0.225$ was used since the amount of data used for estimating the covariate effects for each event type was significantly less. The censoring parameter $c$ was set to be 5. The weighted average estimate cannot always be calculated since the data for each event type could be too sparse to permit an estimate for each type. We only report those estimates which exist. From Table 2, we can see that the weighted average estimator has smaller RASE than that for the local

TABLE 2
*Comparison of the local pseudo-partial likelihood estimator (P) and the weighted average estimator (W)*

|  |  | $\hat{\beta}_1(\cdot)$ | | | | $\hat{\beta}_2(\cdot)$ | | | |
|---|---|---|---|---|---|---|---|---|---|
|  | EST | Abias | SD | SE | RASE | Abias | SD | SE | RASE |
| $\theta = 0.25$ | P | 0.0693 | 0.1063 | 0.1045 | 0.1269 | 0.0198 | 0.1068 | 0.1036 | 0.1086 |
|  | W | 0.0767 | 0.0928 | 0.1222 | 0.1204 | 0.0494 | 0.0898 | 0.1047 | 0.1025 |
| $\theta = 4.00$ | P | 0.0653 | 0.1065 | 0.1029 | 0.0156 | 0.0215 | 0.1041 | 0.1013 | 0.1063 |
|  | W | 0.0750 | 0.0915 | 0.1437 | 0.0140 | 0.0477 | 0.0879 | 0.1037 | 0.1000 |

*Note*: Abias is the average absolute bias of the estimator $\hat{\beta}_j$ for $j = 1, 2$ and RASE denotes the square root of average square errors of the estimator $\hat{\beta}_j$.

TABLE 3
*Comparison of the average square errors of the local pseudo-partial likelihood estimator (P) with those of the one-step estimator (OS)*

| $h$ | Estimator | $\theta = 0.25$, $c = 2$ | | | $\theta = 4.0$, $c = 2$ | | |
|-----|-----------|------|--------|--------|------|--------|--------|
|     |           | mean | median | std    | mean | median | std    |
| 0.1 | P         | 0.0879 | 0.0738 | 0.0602 | 0.0710 | 0.0569 | 0.0550 |
|     | OS        | 0.0879 | 0.0740 | 0.0603 | 0.0716 | 0.0573 | 0.0551 |
| 0.2 | P         | 0.0276 | 0.0199 | 0.0239 | 0.0287 | 0.0220 | 0.0273 |
|     | OS        | 0.0276 | 0.0198 | 0.0239 | 0.0287 | 0.0220 | 0.0272 |
| 0.4 | P         | 0.0107 | 0.0084 | 0.0020 | 0.0216 | 0.0140 | 0.0214 |
|     | OS        | 0.0107 | 0.0084 | 0.0020 | 0.0216 | 0.0140 | 0.0214 |
|     |           | $\theta = 0.25$, $c = 5$ | | | $\theta = 4.0$, $c = 5$ | | |
| 0.1 | P         | 0.0350 | 0.0279 | 0.0256 | 0.0584 | 0.0499 | 0.0506 |
|     | OS        | 0.0350 | 0.0278 | 0.0255 | 0.0586 | 0.0500 | 0.0506 |
| 0.2 | P         | 0.0200 | 0.0137 | 0.0182 | 0.0205 | 0.0148 | 0.0199 |
|     | OS        | 0.0200 | 0.0137 | 0.0182 | 0.0204 | 0.0148 | 0.0199 |
| 0.4 | P         | 0.0162 | 0.0113 | 0.0148 | 0.0162 | 0.0117 | 0.0146 |
|     | OS        | 0.0162 | 0.0113 | 0.0148 | 0.0162 | 0.0117 | 0.0146 |

*Note*: "mean," "median" and "std" denote the average, median and sample standard derivation of the average square errors, respectively, based on 300 simulations.

pseudo-partial likelihood estimator in most of the cases when the weighted average estimator can be calculated.

In the second set of simulations, we compare the performance of the one-step estimator (OS) to that of the maximum local pseudo-partial likelihood estimator (P). We use model (11) with somewhat different configurations. In particular, $V$ is now generated from the standard uniform distribution over $[0, 1]$, $Z$ is independently generated from a standard normal distribution, $\theta = 0.25$ and 4 and $(\lambda_{01}^*, \lambda_{02}^*, \lambda_{03}^*) = (0.2, 1.0, 1.5)$. Censoring times $C_{ij}$ are generated from the uniform distribution over $(0, c)$ with $c = 2$ and 5. We take $g(u) = 8u(1 - u)$ and $\beta(u) = \exp(2u - 1)$.

Table 3 presents the summary of the average square errors ($\text{ASE} = \text{RASE}^2$) for the one-step estimator (OS) and the maximum local pseudo-partial likelihood estimator (P) under various realizations. From the table, we can see that the performance of the one-step estimator is very close to that of the maximum local pseudo-partial likelihood estimator. Figure 2, which presents the box plots for the two methods, also confirms this. We have also conducted simulations using the parameters considered in the first set of simulations. The results are similar.

4.2. *Application to Busselton population health surveys.* We illustrate the proposed method by analyzing a data set from the Busselton Population Health Surveys. The Busselton Population Health Surveys are a series of cross-sectional health surveys conducted in the town of Busselton in Western Australia. Every three years from 1966 to 1981, general health information for each adult participant was collected by means of a questionnaire and a clinical visit. Details of the study are described in [12, 19]. Data for several cardiovascular risk factors are available for 2202 persons who make up 619 families. In this analysis, we investigate the effect of cardiovascular risk factors on the risk of death due to cardiovascular disease (CVD) based on these family data. Since the death times of the family members might be correlated due to genetic factors and cohabitation, we are dealing with multivariate failure time data.

The risk factors we considered included age (in years), body mass index (BMI, in kg/m$^2$), serum cholesterol level (chol) and smoking status. Serum cholesterol (in mmol/L) was determined from a blood sample. Participants' smoking statuses were classified into three categories: never smoker, ex-smoker and current smoker. Two indicator variables were created to indicate the three levels of smoking status: smoke1 is coded as 1 for ex-smoker and 0 otherwise; smoke2 is coded as 1 for current smoker and 0 otherwise.

If a person took part in more than one of the Busselton surveys, only the record from the survey at which that person's age was closest to forty-five years is included. Forty-eight percent of the participants were males (gender = 0 for male and 1 for female). The average age in the data analyzed was 41.7 years, ranging from 16.3 to 89.0 years old. The average cholesterol reading was 5.65 mmol/L. The average body mass index was 24.8 kg/m$^2$. The prevalences of the never-smokers, ex-smokers and current smokers were 49%, 17% and 34%, respectively. Of the 619 families, there were 154 families with one event, 28 families with two events and 3 families with more than two events. There were 219 observed events in all.

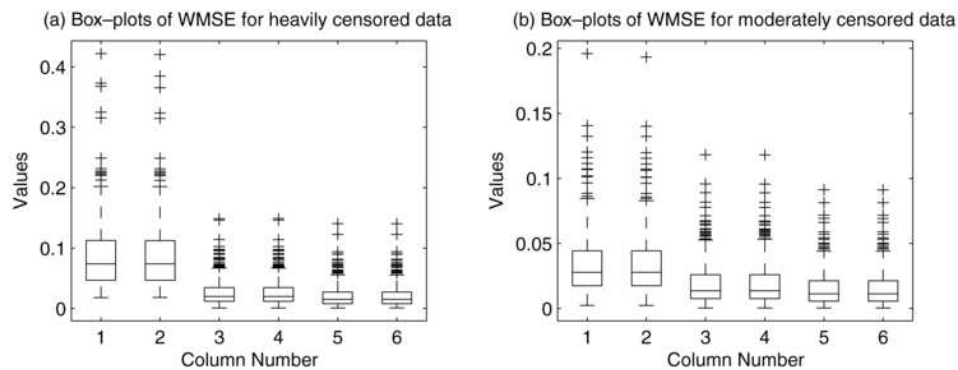

FIG. 2. *Simulation results for the comparison of the maximum local pseudo-partial likelihood estimator (P) with the one-step estimator (OS).* (a) *The box plots are for the distribution of the ASE over the* 300 *replications, using three bandwidths* $h = 0.1, 0.2, 0.4$ *(from left to right). Column numbers* 1, 3 *and* 5 *plot the maximum local pseudo-partial likelihood estimator and column numbers* 2, 4 *and* 6 *plot the one-step estimator (OS) for heavily censored data with* $c = 2$; (b) *the same as* (a) *for moderately censored data,* $c = 5$.

For this analysis, we are interested in investigating how the effect of the risk factors changes with age. We consider the model

$$\lambda_{ij}(t|\mathcal{F}_{ij}) = \lambda_{0j}(t)\exp(\beta_1(\text{age}_{ij}) * \text{gender}_{ij} + \beta_2(\text{age}_{ij}) * \text{bmi}_{ij}$$
$$+ \beta_3(\text{age}_{ij}) * \text{chol}_{ij} + \beta_4(\text{age}_{ij}) * \text{smoke1}_{ij}$$
$$+ \beta_5(\text{age}_{ij}) * \text{smoke2}_{ij} + g(\text{age}_{ij})),$$

where $j = 1$ and 2 denote the parents and the children of the family, respectively, and smoke1 and smoke2 are the indicators for ex-smoker and current smoker, respectively. We take the bandwidth to be $h = 0.15 * (\max(\text{age}) - \min(\text{age})) = 10.905$.

Figure 3 presents the estimates for the varying coefficients as functions of age. From Figure 3(a), we can see that men have a higher risk of dying from CVD than women with the hazard ratio being 1.96 with 95% confidence interval (CI) of $(1.30, 3.03)$ at age fifty. The effect does not seem to change much with age for those older than thirty-five. From Figure 3(b), BMI has little effect on the risk because the coefficient is close to zero over the span of age. From Figure 3(c), higher cholesterol level is associated with higher risk of dying from CVD and the effect of cholesterol increases with age. The hazard ratio for 1 mmol/L change in cholesterol is 1.01 (95% CI: $[0.80, 1.28]$) at age forty and 1.30 (95% CI: $[1.12, 1.53]$) at age sixty-five. From Figures 3(d) and (e), ex-smokers have risk similar to that of the never smokers, while current smokers have a higher risk of dying from CVD compared to the never smokers. The effect of current smoking is higher for younger people with the hazard ratio being 5.60 (95% CI: $[2.19, 14.34]$) at age forty and 1.07 (95% CI: $[0.72, 1.60]$) at age sixty-five.

**5. Concluding remarks.** The local pseudo-partial likelihood is a powerful and a straightforward tool for analyzing multivariate failure time data. The estimator asymptotically follows a normal distribution. Simulation results show that the asymptotic approximation is applicable to finite samples with moderate numbers of clusters.

The weighted average estimator, when it can be calculated, can be a more efficient alternative to the maximum local pseudo-partial likelihood estimator. A disadvantage of the weighted average estimator is that it cannot always be calculated since it involves estimating the covariate coefficient for each failure type and the data for each failure type could be too sparse to permit a reliable estimate.

Use of the one-step estimator is an effective means to reduce the computational burden of an estimator involving iterations. We showed theoretically and empirically that the one-step estimator is an excellent approximation of the fully iterated maximum local pseudo-partial likelihood estimator.

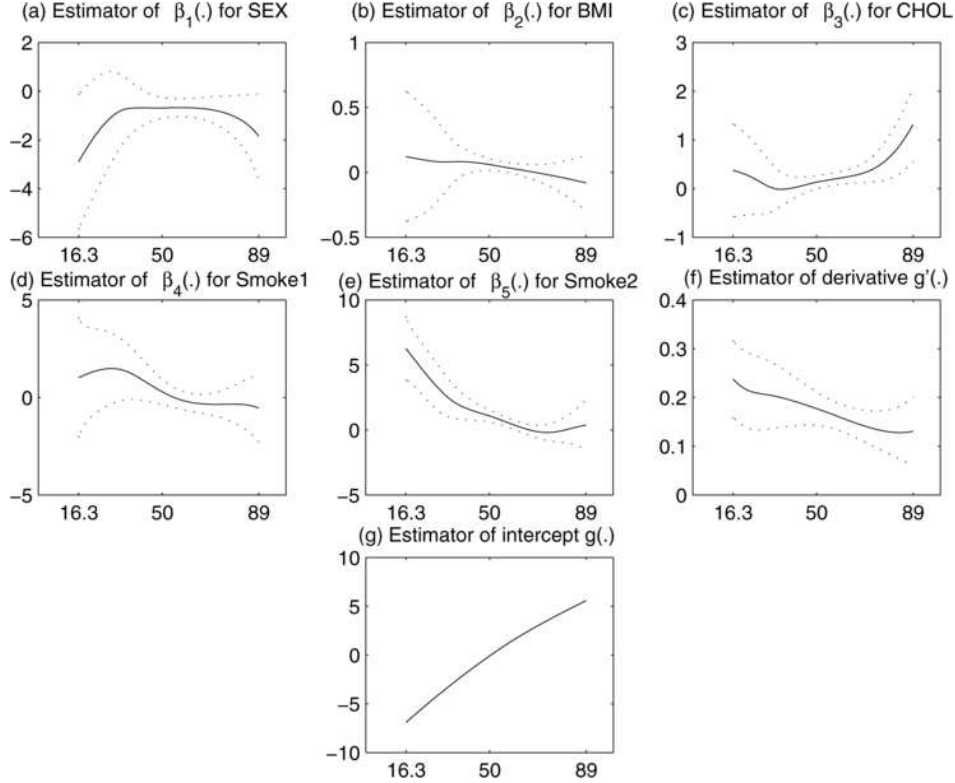

Fig. 3. *Data analysis for Busselton Population Health Surveys study. The marginal hazard rate model is $\lambda_{ij}(t) = \lambda_{0j}(t)\exp(\sum_{k=1}^{5} \beta_k(V_{ij}(t))Z_{ij}(t) + g(V_{ij}(t)))$, where $V(t) = $ age and $\mathbf{Z}^T(t) = $ (Gender, BMI, CHOL, Smoke1, Smoke2), corresponding to the plots (a)–(e), respectively. (f) is the plot of $\hat{g}'(\cdot)$, and (g) is the plot of $\hat{g}(\cdot)$, where smoke1 is coded as 1 for ex-smoker and 0 otherwise, smoke2 is coded as 1 for current smoker and 0 otherwise. The dotted curve is the confidence curve on nominal level $\alpha = 0.05$. In this setting, the chosen bandwidth is $h = 0.15(\max(\mathrm{age}) - \min(\mathrm{age})) = 10.905$. The x-axis is for age.*

Our proposed methods are sensitive to the choice of bandwidth in constructing a local smooth estimation. The one-step estimator and maximum local pseudo-partial likelihood estimator have the same asymptotic distribution and share the same asymptotic bandwidth. We can use the sophisticated bandwidth selection rule proposed in [14] for these estimators.

The methods proposed in this paper can easily be extended to the more general form of multivariate failure time data. More specifically, suppose that there are $n$ clusters that in each cluster there are $K$ correlated individuals and for each individual there are $J$ possible distinct failure types. The marginal hazard function for the $j$th type of failure of the $k$th individual in the $i$th cluster is related to the corresponding covariate vector $\mathbf{Z}_{ijk}(t)$ by

$$\lambda_{ijk}(t, \mathbf{Z}_{ijk}) = \lambda_{0j}(t) \exp\{\beta^T(W_{ijk})\mathbf{Z}_{ijk}(t) + g(W_{ijk})\},$$

where $\lambda_{0j}(t)$, $j = 1, 2, \ldots, J$, are unspecified positive functions and $\beta(\cdot)$ and $g(\cdot)$ are defined as in (1). The maximum local pseudo-partial likelihood for the more general model can be derived similarly and the asymptotic properties can be established with a similar, but more tedious, approach.

**6. Proofs.** Let $(\Omega, \mathcal{F}, P_{(\beta,g,\lambda)})$ be a family of complete probability spaces provided with a history $\mathcal{F} = \{\mathcal{F}_t\}$ for an increasing right-continuous filtration $\mathcal{F}_t \subset \mathcal{F}$. Let $Y_{ij}(t) = I(X_{ij} \geq t)$. We assume that $V_{ij}$ is $\mathcal{F}_{t_{ij}}$-measurable and that $N_{ij}(w)$ and $\mathbf{Z}_{ij}(w)$ are $\mathcal{F}$-adapted. Let $\mathcal{F}_{t,ij} = \sigma\{N_{ij}(w^-), \mathbf{Z}_{ij}(w), V_{ij}, Y_{ij}(w), 0 \leq w \leq t\}$, $i = 1, 2, \ldots, n$, $j = 1, 2, \ldots, J$, and $M_{ij}(t) = N_{ij}(t) - \int_0^t Y_{ij}(w)\lambda_{ij}(w)\,dw$, $i = 1, 2, \ldots, n$. Obviously, $M_{ij}(t)$ is a $\cup_{i=1}^n \mathcal{F}_{t,ij}$ martingale.

To facilitate technical arguments, we will reparameterize the local pseudo-partial likelihood (5) via the transform $\boldsymbol{\zeta} = \mathbf{H}(\boldsymbol{\xi} - \boldsymbol{\xi}_0)$. Hence, the logarithm of the local pseudo-partial likelihood function has the form $\tilde{\ell}_n(\boldsymbol{\zeta}, t) = \ell_n(\mathbf{H}^{-1}\boldsymbol{\zeta} + \boldsymbol{\xi}_0, t)$. By simplification, we have

$$\tilde{\ell}_n(\boldsymbol{\zeta}, t) = \frac{1}{n} \sum_{i=1}^n \sum_{j=1}^J \int_0^t K_h(V_{ij} - v) \\ \times [\boldsymbol{\zeta}^T \mathbf{U}_{ij}^*(w) + \boldsymbol{\xi}_0^T \mathbf{X}_{ij}^*(w) - \log S_{nj0}(w, \boldsymbol{\zeta}, v)]\, dN_{ij}(w),$$

where $\mathbf{U}_{ij}^*(w) = \mathbf{H}^{-1} \mathbf{X}_{ij}^*(w)$ and

$$S_{njk}(w, \boldsymbol{\zeta}, v) = \frac{1}{n} \sum_{i=1}^n K_h(V_{ij} - v) Y_{ij}(w) \exp\{\boldsymbol{\zeta}^T \mathbf{U}_{ij}^*(w) + \boldsymbol{\xi}_0^T \mathbf{X}_{ij}^*(w)\} (\mathbf{U}_{ij}^*(w))^{\otimes k}.$$

Furthermore, for each $w \in [0, \tau]$ and $k = 0, 1, 2$, we rewrite $\tilde{\ell}_n(\boldsymbol{\zeta}) = \tilde{\ell}(\boldsymbol{\zeta}, \tau)$ and define

$$S_{njk}^*(w, \boldsymbol{\theta}, v) = \frac{1}{n} \sum_{i=1}^n K_h(V_{ij} - v) Y_{ij}(w) \exp(\boldsymbol{\beta}^T(V_{ij})\mathbf{Z}_{ij}(w) + g(V_{ij}))(\mathbf{U}_{ij}^*(w))^{\otimes k},$$

where $\boldsymbol{\xi}(\cdot) = (\boldsymbol{\beta}^T(\cdot), \boldsymbol{\beta}'(\cdot)^T, g(\cdot))^T$, $\boldsymbol{\theta}(\cdot) = (\boldsymbol{\beta}^T(\cdot), g(\cdot))^T$, $v \in \Phi_V$.

Recall the notation $\rho(w, \mathbf{z}, v)$ introduced in Section 2.2. Define, for $v \in \Phi_{V,\varepsilon}$,

$$s_{j0}^*(w, \boldsymbol{\theta}, v) = f_j(v) E[\rho(w, \mathbf{Z}_j(w), v) | V_j = v],$$
$$s_{j1}^*(w, \boldsymbol{\theta}, v) = f_j(v) E[\rho(w, \mathbf{Z}_j(w), v)(\mathbf{Z}_j^T(w), 0, 0)^T | V_j = v],$$
$$s_{j2}^*(w, \boldsymbol{\theta}, v) = f_j(v) E[\rho(w, \mathbf{Z}_j(w), v)\Xi(\mathbf{Z}_j, w) | V_j = v],$$

where
$$\Xi(\mathbf{Z}_j, w) = \begin{pmatrix} \mathbf{Z}_j(w)\mathbf{Z}_j^T(w) & \mathbf{0} & \mathbf{0} \\ \mathbf{0} & \mathbf{Z}_j(w)\mathbf{Z}_j^T(w)\mu_2 & \mathbf{Z}_j(w)\mu_2 \\ \mathbf{0} & \mathbf{Z}_j^T(w)\mu_2 & \mu_2 \end{pmatrix}.$$

To facilitate notation, the true functions $\boldsymbol{\theta}_0(u) = (\boldsymbol{\beta}_0^T(u), g_0(u))^T$, $\boldsymbol{\xi}_0(u)$, $\boldsymbol{\zeta}_0 = 0$ and $v$ shall be omitted in $S_{njk}^*(t, \boldsymbol{\theta}, v)$, $S_{njk}(t, \boldsymbol{\zeta}, v)$ and $s_{jk}^*(t, \boldsymbol{\theta}, v)$, $s_{jk}(t, \boldsymbol{\zeta}, v)$ whenever there is no ambiguity, for example,

$$S_{njk}^*(t) = S_{njk}^*(t, v) = S_{njk}^*(t, \boldsymbol{\theta}_0, v), \qquad s_{jk}^*(t) = s_{jk}^*(t, v) = s_{jk}^*(t, \boldsymbol{\theta}_0, v),$$
$$S_{njk}(t) = S_{njk}(t, v) = S_{njk}(t, \mathbf{0}, v), \qquad s_{jk}(t) = s_{jk}(t, v) = s_{jk}(t, \mathbf{0}, v),$$
$$S_{njk}(t, \boldsymbol{\zeta}) = S_{njk}(t, \boldsymbol{\zeta}, v), \qquad\qquad s_{jk}(t, \boldsymbol{\zeta}) = s_{jk}(t, \boldsymbol{\zeta}, v).$$

We will need the following two lemmas in the proofs of the theorems. Let

$$C_{nj}(t) = n^{-1} \sum_{i=1}^{n} Y_{ij}(t)\psi(V_{ij}, (V_{ij} - v)/h, \mathbf{Z}_{ij}(t)) K_h(V_{ij} - v)$$

for a function $\psi(\cdot, \cdot, \cdot)$.

LEMMA 1. *Assume that conditions* (i) *and* (iv) *hold. Assume that* $\psi(\cdot, \cdot, \cdot)$ *is continuous for its three arguments and that* $E(\psi(V_j, w, \mathbf{Z}_j(t))|V_j = v)$ *is continuous at the point* $v$ *for each* $j$ *and* $w$. *If* $h \to 0$ *in such a way that* $nh/\log n \to \infty$, *then*

$$\sup_{0 \le t \le \tau} \sum_{j=1}^{J} |C_{nj}(t) - C_j(t)| \xrightarrow{P} 0,$$

*where* $C_j(t) = f_j(u) \int E(Y(t)\psi(v, w, \mathbf{Z}_j(t))|V_j = v) K(w) \, dw$, $f_j(u)$ *being the density function of* $V$. *Under conditions* (viii)–(x), *we have*

$$\sup_{0 \le t \le \tau} \sup_{v \in B} \sum_{j=1}^{J} |C_{nj}(t, v) - C_j(t, v)| \xrightarrow{P} 0,$$

*where* $B$ *is a compact set satisfying* $\inf_{u \in B} f_j(u) > 0$.

PROOF. By the assumption on $h$, it is easy to show that for every $t \in [0, \tau]$,

(14) $$|C_{nj}(t) - C_j(t)| \xrightarrow{P} 0.$$

Now, we divide $[0, \tau]$ into $M$ subintervals $[t_{i-1}, t_i]$ with a given length not exceeding $\delta$. Note that $\delta$ does not depend on $n$. Then

(15) $$\max_{1 \le i \le M} |C_{nj}(t_i) - C_j(t_i)| \xrightarrow{P} 0.$$

Thus, we need only deal with the term

$$\max_{1\leq i\leq M} \sup_{|t-t_{i-1}|<\delta} |C_{nj}(t) - C_j(t) - \{C_{nj}(t_{i-1}) - C_j(t_{i-1})\}|. \tag{16}$$

By decomposing $\psi$ into a positive part and negative part, we can decompose $C_{nj}(t)$ into $C_{nj}^+(t)$ and $C_{nj}^-(t)$. Hence, we need only show that as $n \to \infty$ and $\delta \to 0$,

$$\max_{1\leq i\leq M} \sup_{|t-t_{i-1}|<\delta} |C_{nj}^+(t) - C_{nj}^+(t_{i-1})|$$
$$+ \max_{1\leq i\leq M} \sup_{|t-t_{i-1}|<\delta} |C_j^+(t) - C_j^+(t_{i-1})| \xrightarrow{P} 0 \tag{17}$$

and a similar result for $C_{nj}^-(t)$ and $C_j^-(t)$.

We now focus on the first term of (17), which is bounded by $J_1 + J_2$, where

$$J_1 = \max_{1\leq i\leq M} \sup_{|t-t_{i-1}|\leq\delta} \left| n^{-1} \sum_{l=1}^n Y_{lj}(t) K_h(V_{lj} - v) \right.$$
$$\times \{\psi^+(V_{lj}, (V_{lj} - v)/h, \mathbf{Z}_{lj}(t))$$
$$\left. - \psi^+(V_{lj}, (V_{lj} - v)/h, \mathbf{Z}_{lj}(t_{i-1}))\} \right|$$

and

$$J_2 = \max_{1\leq i\leq M} \sup_{|t-t_{i-1}|\leq\delta} \left| n^{-1} \sum_{l=1}^n \{Y_{lj}(t) - Y_{lj}(t_{i-1})\} \right.$$
$$\left. \times \psi^+(V_{lj}, (V_{lj} - v)/h, \mathbf{Z}_{lj}(t_{i-1})) K_h(V_{lj} - v) \right|.$$

Note that $\mathbf{Z}_{lj}(t)$, $l = 1, 2, \ldots, n$, is continuous on $[0, \tau]$. Thus, $J_1$ is bounded by

$$\max_{1\leq l\leq n} \max_{1\leq i\leq M} \sup_{|t-t_{i-1}|\leq\delta} |\psi^+(V_{lj}, (V_{lj} - v)/h, \mathbf{Z}_{lj}(t))$$
$$- \psi^+(V_{lj}, (V_{lj} - v)/h, \mathbf{Z}_{lj}(t_{i-1}))| n^{-1} \sum_{l=1}^n K_h(V_{lj} - v),$$

which tends to zero in probability. Since $Y_{lj}(t)$ is a decreasing function of $t$, we have, for any $\varepsilon > 0$, that the probability $P(J_2 > \varepsilon)$ is bounded by

$$MP\left(n^{-1}\left|\sum_{l=1}^n I(t_{i-1} < X_{lj} < t_i)\right.\right.$$
$$\left.\left.\times \psi^+(V_{lj}, (V_{lj} - v)/h, \mathbf{Z}_{lj}(t_{i-1})) K_h(V_{lj} - v)\right| > \varepsilon\right).$$

It is easy to show that

$$n^{-1}\sum_{l=1}^{n} I(t_{i-1} < X_{lj} < t_i)\psi^+(V_{lj}, (V_{lj} - v)/h, \mathbf{Z}_{lj}(t_{i-1}))K_h(V_{lj} - v)$$
$$\xrightarrow{P} \int E\{I(t_{i-1} < X_{lj} < t_i)\psi^+(v, w, \mathbf{Z}_{lj}(t_{i-1}))|V_{lj} = v\}K(w)\,dw.$$

Note that by the Cauchy–Schwarz inequality,

$$E(I(t_{i-1} < X_j < t_i)\psi^+(v, w, \mathbf{Z}_j(t_{i-1}))|V_j = v)$$
$$= |P(t_{i-1}|V_j = v) - P(t_i|V_j = v)|^{1/2} E^{1/2}(\psi^{+2}(v, w, \mathbf{Z}_j(t_{i-1}))|V = v)$$
$$< M^* \delta^{1/2},$$

since $|t_i - t_{i-1}| \leq \delta$, where $M^*$ is some constant. Hence, $J_2 \to 0$ as first $n \to \infty$ and then $\delta \to 0$.

The second term of (17) is bounded by

(18)
$$\max_{1\leq i\leq M} \sup_{|t-t_{i-1}|\leq \delta} f_j(v) \left| \int E\{Y_j(t)(\psi^+(v, w, \mathbf{Z}_j(t)) - \psi^+(v, w, \mathbf{Z}_j(t_{i-1}))|V_j = v\}K(w)\,dw \right|$$
$$+ \max_{1\leq i\leq M} f_j(v) \sup_{|t-t_{i-1}|\leq \delta} \left| \int E\{I(t_{i-1} < X_j < t_i) \times \psi^+(v, w, \mathbf{Z}_j(t_{i-1}))|V_j = v\}K(w)\,dw \right|,$$

which tends to zero as $\delta \to 0$. This implies that (17) holds and hence completes the proof of Lemma 1. $\square$

LEMMA 2. *Under conditions* (i)–(vi), *we have for* $k = 0, 1, 2$,

$$S^*_{njk}(w) = s^*_{jk}(w) + o_p(1),$$

*uniformly for* $w \in (0, \tau]$, *where* $s^*_{jk}(w) = s^*_{jk}(w, \boldsymbol{\theta}_0, v)$ *and, in addition, under conditions* (viii)–(x), *we have*

$$\sup_{w\in(0,\tau], v\in \Phi_V} \|S^*_{njk}(w, v) - s^*_{jk}(w, v)\| = o_p(1).$$

*Furthermore, we have*

$$\sup_{w\in(0,\tau]} \|S_{njk}(w, \boldsymbol{\zeta}) - s_{jk}(w, \boldsymbol{\zeta})\| = o_p(1).$$

Lemma 2 can be easily proven by Lemma 1.

PROOF OF THEOREM 1. We first show that $\hat{\boldsymbol{\zeta}} \to 0$ in probability, where $\hat{\boldsymbol{\zeta}} = \mathbf{H}(\hat{\boldsymbol{\xi}} - \boldsymbol{\xi}_0)$, $\hat{\boldsymbol{\xi}}$ being the maximum local pseudo-partial likelihood estimator of $\boldsymbol{\xi}_0$. Let

$$X_{nj}(t, \boldsymbol{\zeta}) = \frac{1}{n} \sum_{i=1}^{n} \int_0^t K_h(V_{ij} - v) \left[ \boldsymbol{\zeta}^T \mathbf{U}_{ij}^*(w) - \log \frac{S_{nj0}(w, \boldsymbol{\zeta})}{S_{nj0}(w, 0)} \right] dM_{ij}(w).$$

Then it is easy to show that

$$(19) \qquad \tilde{\ell}_n(t, \boldsymbol{\zeta}) - \tilde{\ell}_n(t, \mathbf{0}) = \sum_{j=1}^{J} X_{nj}(t, \boldsymbol{\zeta}) + Y_n(t, \boldsymbol{\zeta}),$$

where

$$Y_n(t, \boldsymbol{\zeta}) = \sum_{j=1}^{J} \int_0^t \left[ (S_{nj1}^*(w))^T \boldsymbol{\zeta} - \log \frac{S_{nj0}(w, \boldsymbol{\zeta})}{S_{nj0}(w, 0)} S_{nj0}^*(w) \right] \lambda_{0j}(w) \, dw.$$

By Lemma 2, we obtain that

$$Y_n(t, \boldsymbol{\zeta}) = \sum_{j=1}^{J} \int_0^t \left[ (s_{j1}^*(w))^T \boldsymbol{\zeta} - \log \frac{s_{j0}(w, \boldsymbol{\zeta})}{s_{j0}(w, 0)} s_{j0}^*(w) \right] \lambda_{0j}(w) \, dw + o_P(1)$$

$$\equiv Y(t, \boldsymbol{\zeta}) + o_P(1).$$

By an argument similar to that in [1], it can be shown that each term in the sum of the asymptotic representation of $Y_n(t, \boldsymbol{\zeta})$ is a strictly concave function in $\boldsymbol{\zeta}$ and that it has the maximum value at $\boldsymbol{\zeta} = 0$. The first term in (19) is a sum of $J$ local square integrable martingales with the square variation process being

$$\langle X_{nj}, X_{nj} \rangle(t) = \frac{1}{n^2} \sum_{i=1}^{n} \int_0^t K_h^2(V_{ij} - v) \left[ \boldsymbol{\zeta}^T \mathbf{U}_{ij}^*(w) - \log \left( \frac{S_{nj0}(w, \boldsymbol{\zeta})}{S_{nj0}(w, 0)} \right) \right]^2$$

$$\times Y_{ij}(w) \exp(\boldsymbol{\beta}_0(V_{ij})^T \mathbf{Z}_{ij}(w) + g_0(V_{ij})) \lambda_{0j}(w) \, dw.$$

It follows from Lemmas 1 and 2 that

$$EX_{nj}^2(t, \boldsymbol{\zeta}) = E \langle X_{nj}, X_{nj} \rangle(t) = O((nh)^{-1}) \to 0, \qquad 0 < t \le \tau.$$

This implies that $X_{nj}(t, \boldsymbol{\zeta}) \to 0$ in probability for $1 \le j \le J$. Hence, we obtain that

$$\tilde{\ell}_n(t, \boldsymbol{\zeta}) - \tilde{\ell}_n(t, \mathbf{0}) = Y(t, \boldsymbol{\zeta}) + o_P(1).$$

We can easily show that $\hat{\boldsymbol{\zeta}}$ maximizes the strictly concave function $\tilde{\ell}_n(t, \boldsymbol{\zeta}) - \tilde{\ell}_n(t, \mathbf{0})$. By Lemma A.1 of [8], it follows that $\hat{\boldsymbol{\zeta}} \to 0$ in probability.

We now prove the second result of Theorem 1. By the same argument as above, we can prove from Lemma 1 that

$$\sup_{t\in[0,\tau]} \sup_{\boldsymbol{\xi}_0 \in \mathbf{C}^*} \sup_{v \in \Phi_V} |\tilde{\ell}_n(t,\boldsymbol{\zeta}) - \tilde{\ell}_n(t,0) - Y(t,\boldsymbol{\zeta})| \longrightarrow 0$$

in probability, where $\boldsymbol{\zeta} = \mathbf{H}(\boldsymbol{\xi} - \boldsymbol{\xi}_0)$ and $\mathbf{C}^*$ is a convex and compact subset of $R^{2p+1}$. Therefore, it follows from Lemma A.1 of [8] that $\sup_{v \in \Phi_V} |\hat{\boldsymbol{\zeta}}| \longrightarrow 0$ in probability. Hence, the proof of Theorem 1 is complete. $\square$

PROOF OF THEOREMS 2 AND 3. Note that we have proved in Theorem 1 that $\mathbf{H}(\hat{\boldsymbol{\xi}}(v) - \boldsymbol{\xi}_0(v)) \to 0$ in probability. This result is very useful for proving Theorem 2. We divide the proofs into the following three steps:

(a) *The asymptotic normality of $\tilde{\ell}'_n(\mathbf{0})$.* The logarithm of the local pseudo-partial likelihood function can be written as

$$\tilde{\ell}'_n(0) = \frac{1}{n} \sum_{i=1}^n \sum_{j=1}^J \int_0^\tau K_h(V_{ij} - v) \left[ \mathbf{U}^*_{ij}(w) - \frac{S_{nj1}(w,v)}{S_{nj0}(w,v)} \right] dM_{ij}(w)$$

$$+ \frac{1}{n} \sum_{i=1}^n \sum_{j=1}^J \int_0^\tau K_h(V_{ij} - v) \left[ \mathbf{U}^*_{ij}(w) - \frac{S_{nj1}(w,v)}{S_{nj0}(w,v)} \right]$$

$$\times \exp(\boldsymbol{\beta}_0(V_{ij})^T \mathbf{Z}_{ij}(w) + g_0(V_{ij})) Y_{ij}(w) \lambda_{0j}(w) \, dw$$

$$\equiv I_1(\tau, 0) + I_2(\tau, 0).$$

We first deal with $I_2(\tau, 0)$. Noting that

$$I_2(\tau,0) = \frac{1}{n} \sum_{i=1}^n \sum_{j=1}^J \int_0^\tau \left( \mathbf{U}^*_{ij}(w) - \frac{S_{nj1}(w)}{S_{nj0}(w)} \right)$$

$$\times [\exp\{\boldsymbol{\beta}_0(V_{ij})^T \mathbf{Z}_{ij}(w) + g_0(V_{ij})\}$$

$$- \exp(\boldsymbol{\xi}_0^T \mathbf{X}^*_{ij} + g_0(v))]$$

$$\times K_h(V_{ij} - v) Y_{ij}(w) \lambda_{0j}(w) \, dw,$$

it follows from a Taylor expansion and Lemma 1 that

$$I_2(\tau,0) = \frac{1}{2n} \sum_{i=1}^n \sum_{j=1}^J \int_0^\tau \left[ \mathbf{U}^*_{ij}(w) - \frac{s^*_{j1}(w)}{s^*_{j0}(w)} \right] Y_{ij}(w) \exp(\boldsymbol{\xi}_0^T \mathbf{X}^*_{ij} + g_0(v))$$

$$\times [\boldsymbol{\beta}''_0(v)^T \mathbf{Z}_{ij}(w) + g''_0(v)](V_{ij} - v)^2 K_h(V_{ij} - v)$$

$$\times \lambda_{0j}(w) \, dw (1 + O_P(h))$$

$$= \frac{1}{2} h^2 \sum_{j=1}^J f_j(v) \int_0^\tau E\left\{ \left[ \begin{pmatrix} \mathbf{Z}_j(w)\mu_2 \\ \mathbf{Z}_j(w)\mu_3 \\ \mu_3 \end{pmatrix} - \frac{s^*_{j1}(w)\mu_2}{s^*_{j0}(w)} \right] \rho(w, \mathbf{Z}_j(w), v) \right.$$

$$\times [\beta_0''(v)^T \mathbf{Z}_j(w) + g_0''(v)] \Big| V_j = v \bigg\}$$
$$\times \lambda_{0j}(w)\, dw(1 + O_P(h)),$$

where $s_{jk}^*(w) = s_{jk}^*(w, \boldsymbol{\theta}_0, v)$ for $k = 0, 1, 2$. Since $K(\cdot)$ is a symmetric function, by simple calculation we have

(20) $\quad I_2(\tau, 0) = \mathbf{B}_n(\tau, v) = \frac{1}{2} h^2 \nu_2 [(\boldsymbol{\Gamma}^{-1} \boldsymbol{\beta}_0''(v))^T, \mathbf{0}^T, 0]^T (1 + O_P(h)).$

We now consider $I_1(\tau, 0)$. Let $\mathbf{B}_{nij}(\tau) = \int_0^\tau K_h(V_{ij} - v)[\mathbf{U}_{ij}^*(w) - \frac{s_{j1}(w, \boldsymbol{\zeta}, v)}{s_{j0}(w, \boldsymbol{\zeta}, v)}] dM_{ij}(w)$. By conditions (vi)–(x), Lemma 2, Lemma A.1 of [25] and some tedious and routine calculation, we can prove that

$$\frac{1}{n} \sum_{i=1}^n \int_0^\tau K_h(V_{ij} - v) \bigg[ \frac{S_{nj1}(w, v)}{S_{nj0}(w, v)} - \frac{s_{j1}(w, v)}{s_{j0}(w, v)} \bigg] dM_{ij}(w) = O_P((nh)^{-1/2}).$$

Hence, it follows that

$$I_1(\tau, 0) = \frac{1}{n} \sum_{i=1}^n \sum_{j=1}^J \mathbf{B}_{nij}(\tau) + o_P(1).$$

Note that $\sqrt{nh} I_1(\tau, 0)$ is a sum of i.i.d. random vectors $\sum_{j=1}^J \mathbf{B}_{nij}(\tau)$ with zero mean and finite variance. The desired asymptotic normality follows from the multivariate central limit theorem by using the Lyapunov condition. It can be shown that the asymptotic variance is

(21)
$$\mathbf{\Pi} = \lim_{n \to \infty} Eh \bigg( \sum_{i=1}^J \mathbf{B}_{nij}(\tau) \bigg)^{\otimes 2}$$
$$= \sum_{j=1}^J \lim_{n \to \infty} Eh \mathbf{B}_{n1j}(\tau)^{\otimes 2} + \sum_{l=1}^J \sum_{j=1, j \neq l}^J \lim_{n \to \infty} Eh \mathbf{B}_{n1j}(\tau) \mathbf{B}_{n1l}(\tau)^T.$$

Note that $\sum_{i=1}^n \mathbf{B}_{nij}(t)$ is a local square-integrable martingale with respect to the filtration $\bigcup_{i=1}^n \mathcal{F}_{t,ij} = \sigma\{N_{ij}(w^-), \mathbf{Z}_{ij}(w), V_{ij}, Y_{ij}(w), 0 \leq w \leq t, i = 1, 2, \ldots, n\}$. Hence, it can be shown that the first term of (21) converges to $\mathbf{D}$. By the Cauchy–Schwarz inequality, we can easily see that $\lim_{n \to \infty} Eh \mathbf{B}_{n1j}(\tau) \times \mathbf{B}_{n1l}(\tau)^T$ exists. Write $\mathbf{\Pi}_{jl}(\tau, v) = \lim_{n \to \infty} Eh \mathbf{B}_{n1j}(\tau) \mathbf{B}_{n1l}(\tau)^T$. Hence, we can prove that the second term of (21) converges to $\mathbf{\Pi}_0(\tau, v) = \sum_{l \neq j} \mathbf{\Pi}_{lj}(\tau, v)$ for the limit matrix $\mathbf{\Pi}_{lj}(\tau, v)$. The proof of Theorem 2 is then completed by using the asymptotic results for $I_1(\tau, 0)$ and $I_2(\tau, 0)$.

(b) *Convergence of the Hessian matrix*. We shall show that the second derivative of the logarithm of the local pseudo-partial likelihood function converges to a finite constant matrix. We have shown in Theorem 1 that $\hat{\boldsymbol{\zeta}} \to 0$ in probability. Hence, by the mean value theorem, we have

(22) $\quad\quad\quad\quad\quad\quad\quad \tilde{\ell}_n''(\hat{\boldsymbol{\zeta}}) = \tilde{\ell}_n''(\mathbf{0}) + o_P(1).$

Since $s_{jk}^*(w) = s_{jk}(w)\exp(g_0(v))$, $k = 0, 1, 2$, from Lemma 2, we can obtain

$$\tilde{\ell}_n''(0) = \frac{1}{n}\int_0^\tau \sum_{i=1}^n \sum_{j=1}^J K_h(V_{ij} - v)\frac{s_{j2}^*(w)s_{j0}^*(w) - s_{j1}^*(w)(s_{j1}^*(w))^T}{(s_{j0}^*(w))^2} dN_{ij}(w)$$
$$+ o_P(1).$$

Write $F_u(w) = P(X \leq w, \Delta = 1|V_j = u)$ and denote its corresponding conditional empirical distribution $\tilde{F}_{nj}(w) = \frac{1}{n}\sum_{i=1}^n K_h(V_{ij} - v)I(X_{ij} \leq w, \Delta_{ij} = 1)$. By means of the conventional argument used in kernel smoothing, together with empirical process theory, it can be shown that

$$\text{(23)} \quad \tilde{\ell}_n''(0) = \sum_{j=1}^J \int_0^\tau \frac{s_{j2}^*(w)s_{j0}^*(w) - s_{j1}^*(w)(s_{j1}^*(w))^T}{(s_{j0}^*(w))^2} d\tilde{F}_{nj}(w)$$
$$= -\mathbf{A}(\tau, v) + o_P(1),$$

where $\mathbf{A}(\tau, v) = \sum_{j=1}^J \int_0^\tau \frac{s_{j2}^*(w)s_{j0}^*(w) - s_{j1}^*(w)(s_{j1}^*(w))^T}{(s_{j0}^*(w))^2} dF_u(w)$. It is easy to show that $\mathbf{A}(\tau, v)$ is positive definite by condition (vii).

(c) *Asymptotic normality of* $\hat{\boldsymbol{\xi}}(v)$. Since $\hat{\boldsymbol{\zeta}}$ maximizes $\tilde{\ell}_n(\boldsymbol{\zeta})$, by Taylor expansion around $\mathbf{0}$, we have

$$-\tilde{\ell}_n'(\mathbf{0}) = \tilde{\ell}_n'(\hat{\boldsymbol{\zeta}}) - \tilde{\ell}_n'(\mathbf{0}) = (\tilde{\ell}_n''(\hat{\boldsymbol{\zeta}}^*))^T\hat{\boldsymbol{\zeta}},$$

where $\hat{\boldsymbol{\zeta}}^*$ lies between 0 and $\hat{\boldsymbol{\zeta}}$ (strictly speaking, the intermediate point can depend on the element of $\tilde{\ell}_n'$, but this does not alter the proof). Hence, $\hat{\boldsymbol{\zeta}}^* \to 0$ in probability. It follows from (23) that

$$\hat{\boldsymbol{\zeta}} - \mathbf{A}(\tau, v)^{-1}\mathbf{B}_n(\tau, v) = -(\tilde{\ell}_n''(\hat{\boldsymbol{\zeta}}^*))^{-1}(\tilde{\ell}_n'(\mathbf{0}) - \mathbf{B}_n(\tau, v)) + o_P(1).$$

By Theorem 1, (23) and Slutsky's theorem, we obtain that

$$\sqrt{nh}(\hat{\boldsymbol{\zeta}} - \mathbf{A}(\tau, v)^{-1}\mathbf{B}_n(\tau, v)) \to N(\mathbf{0}, \mathbf{A}^{-1}(\tau, v)\mathbf{\Pi}(\tau, v)\mathbf{A}^{-1}(\tau, v)).$$

We now simplify the matrix $\mathbf{A}(\tau, v)$. By some simple calculation, we have

$$\text{(24)} \quad s_{j2}^*(w) = \begin{pmatrix} \mathbf{a}_{j2}(w) & \mathbf{0} & \mathbf{0} \\ \mathbf{0} & \mathbf{a}_{j2}(w)\mu_2 & \mathbf{a}_{j1}(w)\mu_2 \\ \mathbf{0} & \mathbf{a}_{j1}^T(w)\mu_2 & \mathbf{a}_{j0}(w)\mu_2 \end{pmatrix}.$$

Similarly, we obtain that $(s_{j1}^*(w))^{\otimes 2} = \text{diag}(\mathbf{a}_{j1}(w)\mathbf{a}_{j1}^T(w), \mathbf{0})$. Note that $s_{j0}^*(w) = \mathbf{a}_{j0}(w)$. By some tedious basic calculation, we have $\mathbf{A}(\tau, v) = \text{diag}(\mathbf{\Gamma}^{-1}, \mathbf{Q}_2\mu_2)$. Hence, the asymptotic bias of the estimator $\hat{\boldsymbol{\xi}}(v)$ is $b(\tau, v) = \mathbf{A}^{-1}(\tau, v)\mathbf{B}_n(\tau, v) = h^2\mathbf{e}_p\boldsymbol{\xi}_0''(v)\mu_2/2$ and the asymptotic covariance is

$$\Sigma(\tau, v) = \mathbf{A}^{-1}(\tau, v)\mathbf{\Pi}(\tau, v)(\mathbf{A}^{-1}(\tau, v))^T$$
$$= \text{diag}(\mathbf{\Gamma}, \mathbf{Q}\mu_2^{-2}\nu_2) + \mathbf{A}^{-1}\mathbf{\Pi}_0(\mathbf{A}^{-1})^T.$$

Therefore, we have finished the proof of the asymptotic normality of the maximum local pseudo-partial likelihood function estimator. □

From the proof of Theorems 2 and 3, we have the following result. If $nh^4 \to 0$, then

$$(25) \quad \sup_{u \in \Phi_{V,\varepsilon}} |(nh)^{1/2}\mathbf{H}\{\hat{\boldsymbol{\xi}}(u) - \boldsymbol{\xi}_0(u)\} - \mathbf{A}^{-1}(\tau, \boldsymbol{\xi}_0(u))\ell'_n(\boldsymbol{\xi}_0(u), \tau)| = o_P(h^{-1/2}).$$

PROOF OF THEOREM 4. By arguments similar to those used in proving Lemma 1, it can be shown that

$$(26) \quad \sup_{t \in [0,\tau]} \sup_{\|\boldsymbol{\theta} - \boldsymbol{\theta}_0\| \leq \|\hat{\boldsymbol{\theta}} - \boldsymbol{\theta}_0\|} n^{-1}|\psi_{nj}(t, \boldsymbol{\theta}) - \psi_{nj}(t, \boldsymbol{\theta}_0)| \longrightarrow 0$$

in probability, where

$$\psi_{nj}(t, \boldsymbol{\theta}) = \sum_{i=1}^{n} I(V_{ij} \in \mathcal{V})Y_{ij}(t)\exp\{\boldsymbol{\beta}^T(V_{ij})\mathbf{Z}_{ij}(t) + g(V_{ij})\},$$

$\boldsymbol{\theta}$ equaling $(\boldsymbol{\beta}^T(\cdot), g(\cdot))^T$. By the definition of $\widehat{\Lambda}_{0j}(t)$, we can show that $\widehat{\Lambda}_{0j}(t) - \Lambda_{0j}(t)$ can be represented by a summation of three terms that are functionals of $\psi_{nj}(t, \boldsymbol{\theta}) - \psi_{nj}(t, \boldsymbol{\theta}_0)$ and it follows that these terms are negligible in the sense of probability. Hence, $\widehat{\Lambda}_{0j}(t) \to \Lambda_{0j}(t)$, uniformly on $(0, \tau]$. Therefore, we can prove by the standard argument of kernel estimation that $\hat{\lambda}_{0j}(t) \to \lambda_{0j}(t)$, uniformly on $(0, \tau]$. □

PROOF OF THEOREM 5. By (25), Theorem 2 and an argument similar to that of Theorem 3 in [25], we can prove Theorem 5. □

**Acknowledgments.** The authors thank Dr. Matthew Knuiman and the Busselton Population Medical Research Foundation in Western Australia for providing the data used in the illustration and are grateful to the Associate Editor and referees for helpful comments.

J. CAI  
H. ZHOU  
DEPARTMENT OF BIOSTATISTICS  
UNIVERSITY OF NORTH CAROLINA AT CHAPEL HILL  
CHAPEL HILL, NORTH CAROLINA 27599-7420  
USA  
E-MAIL: cai@bios.unc.edu  
   zhou@bios.unc.edu  

J. FAN  
DEPARTMENT OF OPERATIONS RESEARCH  
   AND FINANCIAL ENGINEERING  
PRINCETON UNIVERSITY  
PRINCETON, NEW JERSEY 08540  
USA  
E-MAIL: jgfan@princeton.edu  

Y. ZHOU  
INSTITUTE OF APPLIED MATHEMATICS AND  
   CENTER OF STATISTICS  
ACADEMY OF MATHEMATICS AND  
   SYSTEMS SCIENCE  
CHINESE ACADEMY OF SCIENCES  
BEIJING 100080  
CHINA  
E-MAIL: yzhou@amss.ac.cn